\documentclass[10pt]{article}
\usepackage{geometry}
\usepackage{graphicx}
\usepackage{xcolor}
\usepackage{amssymb}
\usepackage{amsmath}
\usepackage{mathtools}
\usepackage[utf8]{inputenc}
\usepackage{amsthm}
\usepackage{bm}
\usepackage{xcolor}

\usepackage{amsthm}

\newtheorem{thm}{Theorem}

\newtheorem*{prf*}{Proof}
\newtheorem{exl}{Example}
\newtheorem{defin}{Definition}

\usepackage{MnSymbol}
\usepackage{algorithmic}
\usepackage{algorithm}
\usepackage{marginnote}
\usepackage{cite}
\usepackage{todonotes}
\usepackage{booktabs}
\usepackage{array}
\usepackage[hidelinks]{hyperref}
\hypersetup{
    colorlinks = false,
    linkbordercolor = {white}
}
\makeatletter
\newcommand{\doubletilde}[1]{{%
  \mathpalette\double@tilde{#1}%
}}
\newcommand{\double@tilde}[2]{%
  \sbox\z@{$\m@th#1\tilde{#2}$}%
  \ht\z@=.9\ht\z@
  \tilde{\box\z@}%
}

\makeatother
\usepackage{amsthm}
\newtheorem{rmrk}{Remark}

\usepackage{mathdots}
\DeclareMathOperator*{\esssup}{ess\,sup} 
\DeclareMathOperator*{\essinf}{ess\,inf} 
\providecommand{\keywords}[1]
{
  \small	
  \textbf{\textit{Keywords---}} #1
}

\begin{document}
	\title{A note on eigenvalues and singular values of variable Toeplitz matrices and matrix-sequences, with application to variable two-step BDF approximations to parabolic equations}

	\author{Nikos Barakitis$^\dag$ \\{ \small\texttt{nickbar@aueb.gr}} \and Valerio Loi$^\dag$ \\{\small\texttt{vloi@uninsubria.it}} \and Stefano Serra-Capizzano$^\dag$ \\{\small\texttt{s.serracapizzano@uninsubria.it}} }


\date{%
    $^\dag$Dipartimento di Scienza ed Alta Tecnologia, Universit\`a dell'Insubria\\
    {\small Via Valleggio 11, 22100 Como, Italy}\\[2ex]%
    \today
}

\maketitle
	
	\begin{abstract}
	The use of variable grid BDF methods for parabolic equations leads to structures that are called variable (coefficient) Toeplitz.
Here, we consider a more general class of matrix-sequences and we prove that they belong to the maximal $*$-algebra of generalized locally Toeplitz (GLT) matrix-sequences. Then, we identify the associated GLT symbols and GLT momentary symbols in the general setting and in the specific case, by providing in both cases a spectral and singular value analysis. More specifically, we use the GLT tools in order to study the asymptotic behaviour of the eigenvalues and singular values of the considered BDF matrix-sequences, in connection with the given non-uniform grids. Numerical examples, visualizations, and open problems end the present work.
	\end{abstract}
	

	\maketitle

\section{Introduction}

The discretization of the continuous on increasingly finer grids for more accurate solutions is the staple and persistent problem in numerical analysis research. This is especially true when considering approximation schemes for differential or integral equations. In fact, because of the invariance in terms of displacement of such operators, after a proper discretization, these problems end up with either Toeplitz or varying Toeplitz matrix-sequences of increasing size. Toeplitz matrices are the matrices with constant main diagonals, while in varying Toeplitz structures the diagonals often vary continuously. In other words, asymptotically, to each diagonal a continuous function is associated and the values are uniform samplings of the given function. While the latter is the most frequent case, also piecewise continuous functions or even Riemann integrable functions are allowed \cite{MR3674485,MR3674485_vol2}.

\par
The Toeplitz matrix-sequences and variations have been widely studied in the past decades in numerous works \cite{BS,Ngondiep}, while the widest generalization of the concept is given by the theory of Generalized Locally Toeplitz (GLT) matrix-sequences, developed in recent years in \cite{Tilli1998,Capizzano2003,glt-2} and in subsequent works (see \cite{GLT-momentary,dyadic-IgA,MR3674485,MR3674485_vol2,Tom} and references therein). A complete presentation of the theory can be found in the books \cite{MR3674485,MR3674485_vol2} and in the papers \cite{barba-max,asympt-unif-gridding,acs-topo}.

Here we first consider a general class of variable Toeplitz matrix-sequences and we prove that they belong to the maximal $*$-algebra of GLT matrix-sequences. The given class contains specific structures stemming from variable two-step backward difference formulae (BDF) approximations of parabolic equations as treated in \cite{BDF-Akrivis}.

Then we identify the associated GLT symbols in the general setting and in the specific case, by providing in both cases a spectral and singular value analysis. We use the GLT tools in order to study the asymptotic behaviour of eigenvalues and singular values of the considered BDF matrix-sequences, also in connection with the given non-uniform grids. Furthermore, we consider a study of positive definiteness and extremal eigenvalues of the involved symmetrized matrices, which is reminiscent of the techniques used for matrix-valued linear positive operators (LPOs) as in \cite{dyadic-cris0,dyadic-LocX,dyadic-IgA}. More specifically, both our theoretical analysis and the numerical tests show that the analysis in \cite{BDF-Akrivis} is sharp in the case of equispaced grids corresponding to Toeplitz structures (see (\ref{Sh-conseq-th2.5})), while there is room for a substantial improvement in the more challenging case of variable griddings associated with a GLT setting.  Numerical examples, visualizations, and open problems end the present work.

The present note is organized as follows.

In Section \ref{sec:tools}, the necessary tools for analysing the spectral behavior of the matrix-sequences of interest are presented, with new material concerning the notion of Toeplitz and GLT momentary symbols originally introduced in \cite{T-momentary,GLT-momentary}. Section \ref{sec:main} is devoted to the study of a general class of variable coefficient matrix-sequences. In Section \ref{sec:BDF}, we specialize our analysis to the case of structures approximating parabolic equations via variable two-step BDF methods: the section contains also a study of positive definiteness and extremal eigenvalues of the involved symmetrized matrices, with tools taken from the theory of matrix-valued LPOs, and numerical experiments for visualizing the theoretical results. A final section of concluding remarks ends the present work.

\section{Spectral tools for matrix analysis}\label{sec:tools}

The section is divided into two parts: first, in Subsection \ref{ssec:distr-vs-structures} we give the definition of spectral and singular value distributions, together with emblematic structures such as Toeplitz, diagonal sampling, and zero-distributed matrix-sequences; subsequently, in Subsection \ref{ssec:GLT}, we present the essentials of the GLT theory.

\subsection{Eigenvalue and Singular Value distribution of matrix-sequences}\label{ssec:distr-vs-structures}

Rather than identifying the eigenvalues or singular values of matrix-sequences of interest exactly which is in most of the cases impossible, their behavior is described in connection with a function $f$ as in Definition \ref{eig-sv-distr}. The approach proves to be adequate for designing fast efficient solvers for a wide class of problems.
\begin{defin}\label{eig-sv-distr}
Let $f:D\subset \mathbb{R}^d \to \mathbb{C}$, $d\ge 1$, be a function such that the Lebesque measure of $D$, denoted by $\mu(D)$, is finite and non-zero. Then the matrix-sequence $\{A_n\}_n$ has an eigenvalue distribution described by $f$, if

\begin{equation}\label{eig-distr}
\lim_{n\to \infty} \frac{1}{n} \sum_{j=1}^{n}F(\lambda_j(A_n)) = \frac{1}{\mu(D)}\int_DF(f(x))dx, \quad \forall F \in C_c(\mathbb{C}),
\end{equation}
where  $C_c(\mathbb{C})$ is the set of all continuous functions with bounded support defined on $\mathbb{C}$. The function $f$ is referred to as the eigenvalue or spectral symbol of the matrix-sequence $\{A_n\}_n$ and we write $\{A_n\}_n \sim_{\lambda}f$. \newline

The matrix-sequence $\{A_n\}_n$ has a singular value distribution described by $f$, if

\begin{equation}\label{sv-distr}
 \lim_{n\to \infty} \frac{1}{n} \sum_{j=1}^{n}F(\sigma_j(A_n)) = \frac{1}{\mu(D)}\int_DF(|f(x)|)dx, \quad \forall F \in C_c(\mathbb{R}),
\end{equation}
where $C_c(\mathbb{R})$ is the set of all continuous with bounded support functions defined on $\mathbb{R}$. In this case the function $f$ is referred to as the singular value symbol of the matrix-sequence $\{A_n\}_n$  and we write $\{A_n\}_n \sim_{\sigma}f$. \newline
If $f$ is both the eigenvalue and singular value symbol of $\{A_n\}_n$ we may write $\{A_n\}_n \sim_{\lambda,\sigma}f$.

\end{defin}

\par
Intuitively speaking, if $d\ge 1$ and the matrix-sequence $\{A_n\}_n$ has an eigenvalue distribution described by $f:D\subset \mathbb{R}^d \to \mathbb{C}$, then under the condition that $f$ is continuous almost everywhere (a.e.) a proper rearrangement of the  eigenvalues of $A_n$ is close to a sampling of $f$ over an equispaced grid on $D$, for $n$ large enough (see e.g. \cite{asympt-unif-gridding} and references there reported). This definition allows eigenvalues to be out of the range of $f$; however, the total number of such eigenvalues is at most $o(n)$.
Furthermore, the case where $f$ is equal to zero a.e in the limit relation (\ref{sv-distr}) identifies the case of zero-distributed matrix-sequences, which represent one of the three building blocks of the GLT matrix-sequences.
	
\subsubsection{Toeplitz matrix-sequences}\label{ssec:toe}

As mentioned above Toeplitz matrices have constant main diagonals. In our setting, each Toeplitz matrix-sequence is associated with a Lebesgue integrable function $f(\theta)\sim \sum_{k=-\infty}^{\infty}a_{k}e^{\mathbf{i}k\theta}$ over $[-\pi,\pi]$, the right part of the equality being the Fourier series of the function, defined on $[-\pi,\pi]$ and periodically extended on the whole real line. The Fourier coefficients $a_k$ of $f$ that is
\begin{equation}\label{eq:fourier_coefficients}
  a_{k} =\frac{1}{2\pi} \int_{[-\pi,\pi]} f({\theta})e^{-\hat i\, k\theta}\, d\theta \in \mathbb{C}, \ \ \ {k}
  \in \mathbb{Z},\ {\hat i}^2=-1,
\end{equation}
for $k=-n+1,\dots,n-1$, are arranged on the diagonals of the matrix of order $n$ as follows
\begin{equation}\label{toeplitz}
	T_{n}(f)=\begin{bmatrix}{}
	a_0 & a_{-1} & \cdots & a_{-n+2} & a_{-n+1} \\
	a_1 & a_0 & a_{-1}   &  & a_{-n+2} \\
	\vdots & a_1 & a_0 & \ddots & \vdots \\
	a_{n-2} &  & \ddots & \ddots & a_{-1} \\
	a_{n-1} & a_{n-2} &\cdots & a_1 & a_0
	\end{bmatrix}.
\end{equation}
The $k$-th Fourier coefficient appears on the diagonal where  the difference row-column index equals $k$. The function $f$ is called the \emph{generating function} of $T_n(f)$ and of the whole matrix-sequence $\{T_n(f)\}_n$.
\par
In \cite{MR890515} Szeg{\H{o}} first proved that the eigenvalues of the Toeplitz matrix-sequence $\{T_n(f)\}_n$ generated by a real-valued function $f\in L^{\infty}([-\pi,\pi])$ are asymptotically distributed as $f$ in the sense of the first part of Definition \ref{eig-sv-distr}. Since then the result has been extended to include real or complex-valued functions  $f\in L^{1}([-\pi,\pi])$ (see \cite{Tilli_Spectr1998,TYRTYSHNIKOV199815} and references therein), with interesting asymptotics on the extreme eigenvalues in the Hermitian setting \cite{SERRA_1996_ON_THE,SERRA1998109,BOTTCHER1998285}. The generalized Szeg{\H{o}} theorem and other findings on the extremal eigenvalues are resumed in the subsequent theorem.
\begin{thm}\label{th:distri-extre}
Suppose that $f \in L^{1}([-\pi,\pi])$ and $\{T_n(f)\}_n$ is the Toeplitz matrix-sequence generated by $f$. Then
		\[
		\{T_n(f)\}_n \sim_{\sigma} f.
		\]
		If  $f$ is real-valued a.e., then the generated matrices are Hermitian and
		\[
		\{T_n(f)\}_n \sim_{\lambda} f.
		\]
		If  $f$ is real-valued a.e., ordering the eigenvalues in non-increasing order and setting $\lambda_j^{(n)}$ to be the $j$-th eigenvalue of $T_n(f)$, then
        \[
        M=\esssup(f)>\lambda_1^{(n)}\geq\lambda_2^{(n)}\geq\dots\geq\lambda_n^{(n)}>\essinf(f)=m \ \ \ \ \forall \: n,
        \]
        if $m<M$. In the case where $m=M$ everything is trivial since $T_n(f)=T_n(m)=mI_n$, with $I_n$ being the identity matrix of size $n$.
		In addition, for any fixed $j$ independent of $n$, we find
		\[
        \lim_{n\rightarrow\infty}\lambda_{j}^{(n)}=M, \quad  \lim_{n\rightarrow\infty}\lambda_{n-j}^{(n)}=m
		\]
		and the convergence speed to $m$ is governed as $n^{-\alpha}$ if $f-m$ has a finite number of zeros, whose maximal order is $\alpha>0$, while, in perfect analogy, the convergence speed to $M$ is governed as $n^{-\beta}$ if $M-f$ has a finite number of zeros, whose maximal order is $\beta>0$.
	
\end{thm}

\subsubsection{Diagonal sampling matrices and matrix-sequences}\label{ssec:diag-sampling}

Let $D\subset\mathbb{R}$ such that it is reasonable to consider an equispaced grid on $D$ and $\alpha:D\rightarrow\mathbb{C}$ a function defined on $D$.
Indeed for having a regular asymptotic behaviour, $D$ is required to be Peano-Jordan measurable which is equivalent to the Riemann integrability of its chacteristic function, while $\alpha$ has to be Riemann integrable over $D$. A diagonal matrix whose diagonal elements are a sampling of $\alpha$ on an equispaced grid on $D$ is called diagonal sampling matrix. Of course, for a matrix-sequence $\{A_n\}_n$ whose diagonal elements are samplings of a function $\alpha$, we have $\{A_n\}_n \sim_{\lambda,\sigma}\alpha$. We here are interested in such matrices where $\alpha$  is a continuous function a.e. on $[0, \: 1]$, i.e. $\alpha$ is Riemann integrable over $[0,1]$. Hence we define $D_n(\alpha)$ as
\[
D_n(\alpha)={\rm diag}_{i=1,\ldots,n}\alpha\Bigl(\frac{i}{n}\Bigr)=\left[\begin{array}{cccc}\alpha(\frac1n) & & & \\ & \alpha(\frac2n) & & \\ & & \ddots & \\   & & & \alpha(1)\end{array}\right].
\]

\subsubsection{Zero-distributed matrix-sequences}\label{ssec:zero-d}

Here we give two crucial definitions and we shortly discuss a connection between them.

\begin{defin} \label{Zero_Distr_Defin}
We say that a matrix-sequence $\{Z_n\}_n$ is zero-distributed if $\{Z_n\}_n \sim_{\sigma} 0$. That is,
\begin{align*}
\lim_{n\rightarrow \infty} \: \frac{1}{n} \sum_{j=1}^{n}F(\sigma_j(Z_n)) = F(0), \ \ \ \  \forall \: F \: \in \: C_c(\mathbb{R}).
\end{align*}
It can be shown that a matrix-sequence $Z_n$ is zero-distributed if and only if
\begin{align}
\forall \: n, \:\: Z_n= R_n + N_n, \quad \text{where,} \:\:\lim_{n\rightarrow\infty}\frac{{\rm rank}(R_n)}{n}=\lim_ {n\rightarrow\infty}\|N_n\|=0. \label{Zero_Distr_Equiv_Defin}
\end{align}
In addition, if $\{Z_n\}_n$ is zero-distributed, and the matrices of the sequence are Hermitian, then $\{Z_n\}_n \sim_{\lambda} 0$.
\end{defin}

The following definition is the main tool for the original construction of GLT matrix-sequences, whose axioms of interest in our work are reported in the subsequent section.

\begin{defin}\label{Approxim_Class_Sequences_Defin}
Let  $\{A_n\}_n$ be a matrix-sequence and  $\{\{B_{n,m}\}_n\}_m $ be a class of matrix-sequences. We say that $\{\{B_{n,m}\}_n\}_m $ is an  approximating class of sequences (a.c.s.)  for $\{A_n\}_n$, and we write $\{B_{n,m}\}_n \xrightarrow{\text{ a.c.s}} \{A_n\}_n$ if, for every $m$, there exists an $n_m$ such that, for $n\geq n_m$,
\begin{align}
A_n=B_{n,m}+R_{n,m}+N_{n,m}, \quad {\rm rank}(R_{n,m}) \leq c(m)n, \quad \|N_{n,m}\| \leq \omega(m),
\end{align}
where $n_m, \: c(m), \: \omega(m)$ depends only on $m$, and $\lim_{m\rightarrow \infty}c(m)=\lim_{m\rightarrow \infty}\omega(m)=0$. Here, $\| \cdot \|$ denotes the spectral norm or induced Euclidean norm, i.e. the Schatten $p$ norm with $p=\infty$ coinciding with the maximal singular value of its argument.
\end{defin}

Notice that (\ref{Zero_Distr_Equiv_Defin}) implies that any zero-distributed matrix-sequence $\{Z_n\}_n$ is such $\{\{B_{n,m}\}_n\}_m $ with $B_{n,m}$ being all the identically zero matrix, converges to $\{Z_n\}_n$ in the a.c.s. topology according to  Definition \ref{Approxim_Class_Sequences_Defin}; see also \cite{acs-topo} and \cite{barba-max} for a deep use of the notion.

\subsection{GLT matrix-sequences}\label{ssec:GLT}
The three special classes of matrix-sequences described above are the main building blocks of any GLT matrix-sequence. These components, in any algebraic combination, including conjugate transposition and inversion if possible generate the whole GLT class \cite{MR3674485}, which forms a maximal $*$-algebra, isometrically equivalent to the $*$-algebra of measurable symbols; see \cite{barba-max}.

\subsection*{Short description of GLT matrix-sequences}
Instead of presenting the complete definition of the GLT class (see \cite{MR3674485}), which is technical and demanding, we here only give an incomplete description of the class and list the axioms, which prove to be sufficient for studying the spectral and singular value distribution of the matrix-sequences of interest (see \cite{MR3674485}[Chapter 9] for the complete set of axioms characterizing uniquely the GLT $*$-algebra). In brief, the GLT class is constructed as follows:
\par
The three classes of matrix-sequences described above, namely the
\begin{itemize}
\item  Toeplitz matrix-sequences generated by a  function in $L^1([-\pi,\pi])$,
\item  diagonal sampling matrix-sequences generated by a Riemann integrable function,
\item  zero-distributed matrix-sequences,
\end{itemize}
belong to the class. Then, whatever results  under the common algebraic operations between GLT matrix-sequences or can be approximated in the a.c.s. sense by such a class of matrix-sequences, also belongs to the class.

 \par

The basic properties of the class are:
\begin{description}
	
	\item[{\bf GLT1}] Every GLT matrix-sequence is related to a unique function $\kappa:[0,\:1]\times[-\pi,\:\pi]\rightarrow \mathbb{C}$, which is the GLT symbol of the matrix-sequence. The singular values of the matrix-sequence are distributed as the function $\kappa$. If the matrices of the sequence are Hermitian, then the eigenvalues of the matrix-sequence are distributed as the function $\kappa$. We denote the GLT matrix-sequence $\{A_n\}_n$ with GLT symbol $\kappa$ as
	\[
	\{A_n\}_n \sim_{\rm GLT} \kappa.
	\]
	\item[{\bf GLT2}] Every Toeplitz matrix-sequence, with generating function $f \in L^1([-\pi,\pi])$ is a GLT matrix-sequence, with symbol $\kappa(x,\theta) = f(\theta)$.
	\item[{\bf GLT3}] Every diagonal matrix-sequence, whose elements are a uniform sampling of a continuous function a.e. $\alpha:[0,\:1]\rightarrow \mathbb{C}$ is a GLT matrix-sequence with symbol $\kappa(x,\theta) = \alpha(x)$.
	\item[{\bf GLT4}] Every zero-distributed matrix-sequence is a GLT matrix-sequence with symbol $\kappa(x,\theta)=0$.
	\item[{\bf GLT5}] The set of all GLT matrix-sequences is a maximal $*$-algebra. In other words, the GLT class is maximal and it is closed under linear combinations, multiplications, conjugate transpositions, and (pseudo)-inversions, provided that the symbol of the matrix-sequence which is (pseudo)-inverted is zero at a set of zero Lebesgue measure. Therefore, a matrix-sequence obtained by operations among GLT matrix-sequences is GLT with a symbol produced by identical operations among the symbols.
	\item[{\bf GLT6}] $\{A_n\}_n \sim_{\rm GLT} \kappa$, if and only if there exist GLT matrix-sequences $\{B_{n,m}\}_n \sim_{\rm GLT} \kappa_m$ such that $\kappa_m$ converge to $\kappa$ in measure and $\{\{B_{n,m}\}_n\}_m$ is an a.c.s. for $\{A_n\}_n$.
\end{description}

{
\begin{rmrk}
For understanding the reason why it is reasonable to use the GLT $*$-algebra in a discretization process it is enough to bring in mind the components of a differential equation. The discrete analogous of the differential operators applied to the unknown function are the Toeplitz matrices while the discrete analogous of the coefficient functions are the diagonal sampling matrices. Very importantly, with respect to the a.c.s. topology, a zero-distributed matrix-sequence is the maximal deviation allowed for two different matrix-sequences having the same eigenvalue (singular value) distribution.
We observe that both Definition \ref{Zero_Distr_Defin} and Definition \ref{Approxim_Class_Sequences_Defin} clarify the two directions and the corresponding tolerance limits of the possible deviations.
\end{rmrk}
}

Now we introduce the notion of Toeplitz and GLT momentary symbols. These are slight variations of Definition 4 in \cite{T-momentary} and Definition 1.1 in \cite{GLT-momentary}, respectively. The novelty relies in the requirement for the remainder term $\{{R_n}\}$ as in (\ref{rel:normaliz-cond-new1}) and (\ref{rel:normaliz-cond-new2}), which is more precise in the present formulation.
\begin{defin}[Toeplitz momentary symbols]
\label{def:momentarysymbols-new}
Let $\{X_n\}_n$ be a matrix-sequence and assume that there exist matrix-sequences $\{A_n^{(j)}\}_n$, $\{R_n\}_n$, scalar sequences $c_n^{(j)}$, $j=0,\ldots,t$, and Lebesgue integrable functions $f_j$ defined over $[-\pi,\pi]$, $t$ nonnegative integer independent  of $n$, such that
\begin{equation} \label{rel:normaliz-cond-new1}
\left\{\frac{R_n}{c_n^{(t)}}  \right\}_n\sim_{\rm GLT} 0
\end{equation}
is zero-distributed (as in item {\textbf{GLT4}}),
\begin{eqnarray} \nonumber
\left\{ \frac{A_n^{(j)}}{ c_n^{(j)}}\right\}_n & = & T_n(f_j), \\ \nonumber 
c_n^{(0)}=1, & & c_n^{(s)}=o(c_n^{(r)}), \ \ t\ge s>r, \\
\{X_n\}_n & = & \{A_n^{(0)}\}_n + \sum_{j=1}^t \{A_n^{(j)}\}_n + \{R_n\}_n, \nonumber
\end{eqnarray}
$\{A_n^{(j)}\}_n$, $j=0,\ldots,t$, not zero-distributed which imply $f_j$,
$j=0,\ldots,t$, not identically zero almost everywhere. Then, by a slight abuse of notation, the function
\begin{equation}\label{eq:T_momentary_1D}
f_n=f_0+ \sum_{j=1}^t c_n^{(j)} f_j
\end{equation}
is defined as the Toeplitz momentary symbol for $X_n$ and $\{f_n\}$ is the sequence of Toeplitz momentary symbols for the matrix-sequence $\{X_n\}_n$.
\end{defin}

\begin{defin}[GLT momentary symbols]
\label{def:GLT momentarysymbols-new}
Let $\{X_n\}_n$ be a matrix-sequence and assume that there exist matrix-sequences $\{A_n^{(j)}\}_n$,  $\{R_n\}_n$, scalar sequences $c_n^{(j)}$, $j=0,\ldots,t$, and measurable functions $f_j$ defined over $[0,1]\times [-\pi,\pi]$, $t$ nonnegative integer independent  of $n$, such that
\begin{equation} \label{rel:normaliz-cond-new2}
\left\{\frac{R_n}{c_n^{(t)}}  \right\}_n\sim_{\rm GLT} 0
\end{equation}
is zero-distributed (as in item {\textbf{GLT4}}),
\begin{eqnarray} \nonumber
&\left\{ \frac{A_n^{(j)}}{ c_n^{(j)}}\right\}_n \sim_{\textsc{glt}} f_j, \\ \nonumber
&c_n^{(0)}=1, \quad c_n^{(s)}=o(c_n^{(r)}), \quad t\ge s>r, \\ \nonumber 
&\{X_n\}_n  =  \{A_n^{(0)}\}_n + \sum_{j=1}^t \{A_n^{(j)}\}_n + \{R_n\}_n,
\end{eqnarray}
$\{A_n^{(j)}\}_n$, $j=0,\ldots,t$, not zero-distributed which imply $f_j$,
$j=0,\ldots,t$, all not identically zero almost everywhere. Then, with a slight abuse of notation, the function
\begin{equation}\label{eq:GLT_momentary_1D}
f_n=f_0+ \sum_{j=1}^t c_n^{(j)} f_j
\end{equation}
is defined as the GLT momentary symbol for $X_n$ and $\{f_n\}$ is the sequence of GLT momentary symbols for the matrix-sequence $\{X_n\}_n$.
\end{defin}

Both definitions have of course their block multilevel versions, both in the block multilevel Toeplitz setting and in the block multilevel GLT setting (see \cite{barb d-dim} for the relevant theoretical background).

\section{Main results}\label{sec:main}

Taking inspiration from the recent work \cite{BDF-Akrivis}, we consider a class of banded variable-Toeplitz matrix-sequences associated with two parameters $\delta$, $\eta$, and with a Riemann integrable function $\phi$ equipped with a uniform grid-sequence over its definition domain $[0,T]$, $T>0$. For such a class of matix-sequences we give the GLT and distribution analysis in the sense of the singular values and eigenvalues in Theorem \ref{th:gen}.

\begin{thm}\label{th:gen}
Let $r_i=\phi(\frac{iT}{n}), \: for\: 2=1,\dots,n$, let $\phi:[0,\:T]\rightarrow \mathbb{R}$ be a real-valued Riemann integrable over $[0,\:T]$ and $\delta$, $\eta$ fixed constants.
Then the matrix-sequence $\{\mathbb{L}_n\}_n$, with
\begin{align}
\mathbb{L}_n=\mathbb{L}(r_2,\dots,r_n)=\left[
\begin{array}{ccccc}
\frac{1}{1+r_2}                 &                                    &                                &                \\
-\delta\frac{\sqrt{r_3}}{1+r_3} &\frac{1}{1+r_3}                     &                                &                \\
-\eta\frac{\sqrt{r_3r_4}}{1+r_4}&-\delta\frac{\sqrt{r_4}}{1+r_4}     &\frac{1}{1+r4}                  &                \\
                                &                                    &                                &                \\
                                &-\eta\frac{\sqrt{r_{n-1}r_n}}{1+r_n}&-\delta\frac{\sqrt{r_n}}{1+r_n} &\frac{1}{1+r_n} \\
\end{array}
\right]
\end{align}
 is a GLT matrix-sequence with GLT symbol given by
\[
\kappa(x,\theta)=\frac{1}{1+\hat{\phi}(x)}\left[1-\delta \sqrt{\hat{\phi}(x)}e^{\hat i\, \theta} -\eta \hat{\phi}(x) e^{\hat i\, 2\theta}\right],
\]
where $\hat{\phi}(x)=\phi(Tx)$ so that the physical variable $x$ is defined in the interval $[0,1]$ as requested by the GLT axioms.
Furthermore, setting $\mathbb{S}_n=\Re(\mathbb{L}_n)=\frac{1}{2}\left( \mathbb{L}_n + \mathbb{L}_n^*\right)$, we deduce that $\{\mathbb{S}_n\}_n$ is a GLT matrix-sequence with GLT symbol given by
\[
\Re(\kappa(x,\theta))=\frac{1}{1+\hat{\phi}(x)}\left[1-\delta \sqrt{{\hat{\phi}}(x)}\cos(\theta) -\eta {\hat{\phi}}(x) \cos(2\theta)\right].
\]
Finally $\{\mathbb{L}_n\}_n\sim_\sigma \kappa(x,\theta)$ and $\{\mathbb{S}_n\}_n\sim_\lambda \Re(\kappa(x,\theta))$.
\end{thm}
\begin{prf*}
Let
\begin{align*}
&D_n^{(1)}={\rm diag}\left(\frac{1}{1+r_2},\frac{1}{1+r_3},\dots,\frac{1}{1+r_n}\right),\\
&D_n^{(2)}={\rm diag}\left(\frac{\sqrt{r_2}}{1+r_2},\frac{\sqrt{r_3}}{1+r_3},\dots,\frac{\sqrt{r_n}}{1+r_n}\right),\\
&D_n^{(3)}={\rm diag}\left(\frac{\sqrt{r_1r_2}}{1+r_2},\frac{\sqrt{r_2r_3}}{1+r_3},\dots,\frac{\sqrt{r_{n-1}r_n}}{1+r_n}\right).
\end{align*}
Then a direct check shows that
\begin{align}
\mathbb{L}_n= D_n^{(1)}-\delta D_n^{(2)}T_{n-1}(e^{i\theta}) - \eta D_n^{(3)}T_{n-1}(e^{i2\theta}).
\end{align}
A technical difficulty is given by the fact that the matrices $D_n^{(j)}$, $j=1,2,3$, are diagonal sampling matrices, but not on the canonical grid given in Subsection \ref{ssec:diag-sampling}. In fact
\begin{align}
    D_n^{(1)}&=D_{n-1}\left(\frac{1}{1+\hat{\phi}(x)}\right) + E_n^{(1)}, \\
    D_n^{(2)}&=D_{n-1}\left(\frac{\sqrt{\hat{\phi}(x)}}{1+\hat{\phi}(x)})\right)+ E_n^{(2)},\\
    D_n^{(3)}&=D_{n-1}\left(\frac{\hat{\phi}(x)}{1+\hat{\phi}(x)}\right) + E_n^{(3)},
\end{align}
with $\{E_n^{(j)}\}_n$, $j=1,2,3$, being all zero-distributed in the sense of Subsection \ref{ssec:zero-d} so that $\{E_n^{(j)}\}_n\sim_{\rm GLT} 0$, $j=1,2,3$, since the canonical grid-sequence and the actual grid-sequence are equidistributed \cite{tau2}, due to the Riemann integrability of $\hat{\phi}$ induced by the Riemann integrability of ${\phi}$; see \cite{asympt-unif-gridding}.

As a consequence, using {\rm[{\bf GLT3}]} and recalling that $\{E_n^{(1)}\}_n$, $\{E_n^{(2)}\}_n$, $\{E_n^{(3)}\}_n$ are all GLT matrix-sequences with $0$ symbol, the GLT symbols of $\{D_n^{(1)}\}_n$, $\{D_n^{(2)}\}_n$, $\{D_n^{(3)}\}_n$ are $\frac{1}{1+\hat{\phi}(x)}$, $\frac{\sqrt{\hat{\phi}(x)}}{1+\hat{\phi}(x)}$, $\frac{\hat{\phi}(x)}{1+\hat{\phi}(x)}$, respectively, in the light of axioms {\rm[{\bf GLT3}]}, {\rm[{\bf GLT4}]}, {\rm[{\bf GLT5}]}. \\
Finally from {\rm[{\bf GLT2}]}, {\rm[{\bf GLT3}]}, {\rm[{\bf GLT5}]}, we deduce that $\{\mathbb{L}_n\}_n$ is a GLT sequence with symbol
\[
\kappa(x,\theta)=\frac{1}{1+{\hat{\phi}}(x)}\left[1-\delta \sqrt{{\hat{\phi}}(x)}e^{\hat i\, \theta} -\eta {\hat{\phi}}(x) e^{\hat i\, 2\theta}\right].
\]
Taking the real part of $\{\mathbb{L}_n\}_n$, $\Re(\mathbb{L}_n)=\frac{1}{2}\left( \mathbb{L}_n + \mathbb{L}_n^*\right)$, using the fact that the symbol of $ \{\mathbb{L}_n^*\}_n$ is $\overline{\kappa(x,\theta)}$ by axiom ${\rm[{\bf GLT5}]}$, we deduce that the GLT symbol of $\{\Re(\mathbb{L}_n)\}_n$ is
\[
\Re(\kappa(x,\theta))=\frac{1}{2}\left(\kappa(x,\theta)+ \overline{\kappa(x,\theta)}\right)=
\frac{1}{1+\hat{\phi}(x)}\left[1-\delta \sqrt{{\hat{\phi}}(x)}\cos(\theta) -\eta {\hat{\phi}}(x) \cos(2\theta)\right],
\]
again by axiom ${\rm[{\bf GLT5}]}$.

With the latter steps the proof is concluded by invoking axiom {\rm[{\bf GLT1}]}, which implies $\{\mathbb{L}_n\}_n\sim_\sigma \kappa(x,\theta)$ and $\{\mathbb{S}_n\}_n\sim_\lambda \Re(\kappa(x,\theta))$, $\mathbb{S}_n=\Re(\mathbb{L}_n)$ being Hermitian for every $n$.
\end{prf*}

\section{A specific setting}\label{sec:BDF}
	
Here we deal in some detail with the specific setting in paper \cite{BDF-Akrivis} that inspired our work.  Following \cite{BDF-Akrivis}, let $T>0,\:u^{(0)} \:\in\:H$, and consider the initial value problem
\begin{align}
\begin{cases}\label{cont pb}
u'(t)+Au(t) = f(t), \quad 0<t<T, \\
u(0)=u^{(0)},
\end{cases}
\end{align}
in which we look for the solution $u\: \in \: C( (0,\:T];\mathcal{D}(A) ))\bigcap C([0,\:T];H)$ satisfying the conditions in (\ref{cont pb}).
Here $A$ represents a positive definite, selfadjoint linear operator on a Hilbert space $(H,(\cdot,\cdot))$, with domain $\mathcal{D}(A)$ dense in $H$ and $f:[0,T]\rightarrow H$ is the given forcing term.
\par
Methods based on BDFs are popular for stiff differential equations, in particular, for parabolic equations. They are frequently implemented on nonuniform partitions for numerical efficiency: in particular a varying grid is crucial for dealing in a different way with time intervals showing fast variations of the solution and others characterized by slow variations of the solution.
\par
For an integer $n\geq 2$, consider a partition $0=t_0<t_1<\dots<t_n=T$ of the time interval $[0,\:T]$, with time steps $k_i:=t_i-t_{i-1}$, $i=1,\dots,n$. We recursively define a sequence of approximations $u^{(i)}$ to the nodal values $u(t_i)$ of the exact solution by the variable two-step BDF method,
\begin{align}
D_2u^{(i)}+Au^{(i)}=f^{(i)}, \quad i=2,\dots,n,
\end{align}
with $f^{(i)}:=f(t_i)$, assuming that arbitrary starting approximations $u^{(0)}$ and $u^{(1)}$ are given. Here,
\[
D_2v^{(i)}:=\left(1+\frac{k_i}{k_{i-1}}\right)\frac{v^{(i)}-v^{(i-1)}}{k_i}-\frac{k_i}{k_{i-1}}\frac{v^{(i)}-v^{(i-2)}}{k_i+k_{i-1}}
\]

With reference to Theorem \ref{th:gen} and to the work of Akrivis et al. \cite{BDF-Akrivis}, we have $\delta=0.9672$, $\eta= -0.1793$, $t_i=\psi(\frac{iT}{n})$, for $i=1,\dots,n$, $\psi:[0,\:T]\rightarrow[0,\:T]$ increasing and such that $\psi(0)=0,\:\psi(T)=T$. Furthermore, the parameters $r_i$ are defined as
\[
r_i=\frac{t_i-t_{i-1}}{t_{i-1}-t_{i-2}}
\]
for $i=2,\ldots,n$.

if we suppose $\psi$ differentiable then we have
\begin{equation}\label{fata-morgana}
r_i=\frac{t_i-t_{i-1}}{t_{i-1}-t_{i-2}}=\frac{\psi(\frac{iT}{n})-\psi(\frac{(i-1)T}{n})}{\psi(\frac{(i-1)T}{n})-\psi(\frac{(i-2)T}{n})}=\frac{\left(\psi(\frac{iT}{n})-\psi(\frac{(i-1)T}{n})\right)/(T/n)}{\left(\psi(\frac{(i-1)T}{n})-\psi(\frac{(i-2)T}{n}\right)/(T/n)}\rightarrow \frac{\psi'(\frac{iT}{n})}{\psi'(\frac{(i-1)T}{n})}=1
\end{equation}
asymptotically.

However the function $\psi(x)$ is differentiable a.e. because it is monotone. Assuming that $\psi$, as previously discussed, is Riemann integrable, it follows that the derivative of $\psi$ is also continuously differentiable a.e. and $r_i \to 1$ a.e. In this context, as a byproduct of Theorem \ref{th:gen}, the symbol of the sequence $\{\mathbb{L}_n\}_n$ is simply defined as
\[
\kappa(x,\theta)=\kappa(\theta)=\frac{1}{2}-\delta \frac{1}{2}e^{\hat i\, \theta} -\eta \frac{1}{2} e^{\hat i\, 2\theta},
\]
which implies that $\{\mathbb{L}_n-T_{n-1}(\kappa(\theta))\}_n \sim_{{\rm GLT}} 0$, $\{\mathbb{S}_n-T_{n-1}(\Re(\kappa(\theta)))\}_n \sim_{{\rm GLT}} 0$ by axioms {\rm[{\bf GLT2}]}, {\rm[{\bf GLT5}]}.

In Theorem \ref{th:specific}, we give more information on the considered specific setting.

\begin{thm}\label{th:specific}
Let $\mathbb{L}_n=\mathbb{L}(r_2,\dots,r_n)$ and $\mathbb{S}_n=\Re(\mathbb{L}_n)=\frac{1}{2}\left( \mathbb{L}_n + \mathbb{L}_n^*\right)$ be defined as in Theorem \ref{th:gen} and let $r_i$, $i=2,\dots,n$, be as in the beginning of this section with $\delta$, $\eta$ given real constants.
Then $\{\mathbb{L}_n\}_n$ is a GLT matrix-sequence with symbol $\kappa(\theta)=\frac{1}{2}-\delta \frac{1}{2}e^{\hat i\, \theta}
-\eta \frac{1}{2} e^{\hat i\, 2\theta}$ and  $\{\mathbb{S}_n\}_n$ is a GLT matrix-sequence with GLT symbol given by
\[
\Re(\kappa(\theta))=\frac{1}{2}\left[1-\delta \cos(\theta) -\eta \cos(2\theta)\right],
\]
so that $\{\mathbb{L}_n\}_n\sim_\sigma \kappa(\theta)$ and $\{\mathbb{S}_n\}_n\sim_\lambda \Re(\kappa(\theta))$.
Finally, under the assumption that the grid is defined by $\psi$ strictly increasing and twice continuously differentiable, both $\{\mathbb{L}_n\}_n$ and $\{\mathbb{S}_n\}_n$ have a GLT momentary expansion with zero terms given by the standard GLT symbols and first terms having the expression
\[
l(x,\theta)=-\frac{\hat{\psi}''(x)}{\hat{\psi}'(x)}\left[\frac{1}{4}+\frac{\eta}{4}e^{2\hat{i}\theta}\right], \ \ \ \Re(l(x,\theta))=-\frac{\hat{\psi}''(x)}{\hat{\psi}'(x)}\left[\frac{1}{4}+\frac{\eta}{4}\cos(2\theta)\right],
\]
respectively. Here $\hat{\psi}(x)=\psi(xT)$ so that it is defined in $[0,1]$ as required by the GLT theory.
\end{thm}
\begin{prf*}
The fact that $\{\mathbb{L}_n\}_n$ is a GLT matrix-sequence with symbol $\kappa(\theta)=\frac{1}{2}-\delta \frac{1}{2}e^{\hat i\, \theta}
-\eta \frac{1}{2} e^{\hat i\, 2\theta}$ and  $\{\mathbb{S}_n\}_n$ is a GLT matrix-sequence with GLT symbol given by $\Re(\kappa(\theta))$ are obvious consequences of Theorem \ref{th:gen} and (\ref{fata-morgana}), while  $\{\mathbb{L}_n\}_n\sim_\sigma \kappa(\theta)$ and $\{\mathbb{S}_n\}_n\sim_\lambda \Re(\kappa(\theta))$ follow from axiom {\rm[{\bf GLT1}]}.

The second part is essentially a consequence of the use of Taylor expansions and of the notion of GLT momentary symbol (see \cite{GLT-momentary} and the new Definition \ref{def:GLT momentarysymbols-new}). More precisely, we identify an expansion of $\{\mathbb{L}_n\}$ that takes the form
\begin{equation}
\{\mathbb{L}_n\}_n = \{A_n\}_n + \sum_{j=1}^t h_n^{(j)} \{B^{(j)}_n\}_n+ \{R_n\}_n,
\end{equation}
where we find a principal term
\begin{equation*}
\left\{ A \right\}_n \sim_{{\rm GLT}} f_0,
\end{equation*}
``higher order terms"
\begin{equation*}
\left\{ {B^{(j)}_n} \right\}_n \sim_{{\rm GLT}} f_j,
\end{equation*}
with
\begin{align*}
h^{(1)}&=o(1),\\
h^{(i+j)}_n & = o(h^{(i)}_n), \quad j>0, \ \ i+j\le t,
\end{align*}
and $\left\{ {R_n\over h_n^{(t)}} \right\}_n$ being a zero-distributed remainder, according to (\ref{rel:normaliz-cond-new2}).

Then, the momentary symbol sequence is the sequence of functions as follows
\begin{equation}
g_n = f_0 + \sum_{j=1}^t h^{(j)}_n f_j.
\end{equation}

Note that, as a consequence of our definition, $f_0$ represents the GLT symbol of $\{\mathbb{L}_n\}_n$, which means $f_0=\frac{1}{2} - \delta \frac{1}{2} e^{\hat{i} \theta} - \eta \frac{1}{2} e^{\hat{i} 2\theta}=\kappa(\theta)$ and $\lim_{n \to \infty} g_n = f_0$ a.e. In our context, we only need a first-order expansion, and to find such an expansion, we rely extensively on the Taylor theorem on the spatial symbol.
Let
\[
r_i = \frac{\psi \left( \frac{iT}{n} \right) - \psi \left( \frac{(i-1)T}{n} \right)}{ \psi \left( \frac{(i-1)T}{n} \right) - \psi \left( \frac{(i-2)T}{n} \right)}.
\]
By the Taylor theorem, we observe that
\begin{align*}
    \psi \left( \frac{iT}{n} \right) &=  \psi \left( \frac{(i-1)T}{n} \right) + h_n \psi' \left( \frac{(i-1)T}{n} \right) + \frac{h_n^2}{2}\psi^{''} \left( \frac{(i-1)T}{n} \right) + o(h_n^2) \\
    \psi \left( \frac{(i-2)T}{n} \right) &=  \psi \left( \frac{(i-1)T}{n} \right) - h_n \psi' \left( \frac{(i-1)T}{n} \right) + \frac{h_n^2}{2}\psi^{''} \left( \frac{(i-1)T}{n} \right) + o(h_n^2),
\end{align*}
where $h_n=\frac{T}{n}$.
This enables us to express the local approximation of $r_i$ as follows

\begin{equation}
    r_i=\frac{h_n \psi' \left( \frac{(i-1)T}{n} \right) + \frac{h_n^2}{2}\psi^{''} \left( \frac{(i-1)T}{n} \right) + O(h_n^2)}{h_n \psi' \left( \frac{(i-1)T}{n} \right) - \frac{h_n^2}{2}\psi^{''} \left( \frac{(i-1)T}{n} \right) + O(h_n^2)}= 1+ h_n\frac{\psi^{''} \left( \frac{(i-1)T}{n} \right)}{\psi^{'} \left( \frac{(i-1)T}{n} \right)} + o(h_n),
\end{equation}
the last equality is a direct consequence of the classical Taylor expansion $\frac{1+x}{1-x}=1+2x+O(x^2)$. \\

Our goal is to utilize this approximation of $r_i$ to determine the first-order momentary GLT expansion of $\{\mathbb{L}_n\}_n$. We recall that

$$
\mathbb{L}_n = D^{(1)}_n \left( I - \delta D^{(2)}_n T_{n-1}(e^{i \theta}) - \eta D^{(3)}_n T_{n-1}(e^{i 2\theta}) \right),
$$
where
\begin{align*}
D^{(1)}_n &= {\rm diag} \left( \frac{1}{1 + r_2}, \frac{1}{1 + r_3}, \ldots, \frac{1}{1 + r_n} \right), \\
D^{(2)}_n &= {\rm diag} \left( {\sqrt{r_2}}, {\sqrt{r_3}}, \ldots, \sqrt{r_n} \right), \\
D^{(3)}_n &= {\rm diag} \left( \sqrt{r_1 r_2}, \sqrt{r_2 r_3}, \ldots, \sqrt{r_{n-1} r_n} \right).
\end{align*}

Using the previously identified approximation and further expanding the samplings on the diagonal matrices, we derive

\begin{align*}
D^{(1)}_n &= {\rm diag} \left( \frac{1}{2} - \frac{h_n}{4}\frac{\psi^{''} \left( \frac{(i-1)T}{n} \right)}{\psi^{'} \left( \frac{(i-1)T}{n} \right)} \right) + F_{1,n}, \\
D^{(2)}_n &= {\rm diag} \left( 1+\frac{h_n}{2} \frac{\psi^{''} \left( \frac{(i-1)T}{n} \right)}{\psi^{'} \left( \frac{(i-1)T}{n} \right)} \right) +
F_{2,n}, \\
D^{(3)}_n &= {\rm diag}\left( 1+ h_n \frac{\psi^{''} \left( \frac{(i-1)T}{n} \right)}{\psi^{'} \left( \frac{(i-1)T}{n} \right)} \right) + F_{3,n},
\end{align*}
where $\left\{ {F_{i,n}\over h_n} \right\}_n$, $i=1,2,3$, are diagonal zero-distributed error matrix-sequences.

By breaking down diagonal matrices into sums, carrying out algebraic operations, and leaving some zero-distributed remainders, we obtain
\begin{align*}
    \mathbb{L}_n &= D^{(1)}_n \left( I - \delta D^{(2)}_n T_{n-1}(e^{\hat{i} \theta}) - \eta D^{(3)}_n T_{n-1}(e^{\hat{i} 2\theta}) \right) \\
    &=\frac{1}{2}I-\frac{\delta}{2}T_{n-1}\left(e^{\hat{i} \theta}\right)-\frac{\eta}{2}T_{n-1}\left( e^{2 \hat{i} \theta} \right)+ h_n\left[ {\rm diag} \left(-\frac{\psi^{''} \left( \frac{(i-1)T}{n} \right)}{\psi^{'} \left( \frac{(i-1)T}{n} \right)}\right)\left(\frac{1}{4}I+\frac{\eta}{4}T_{n-1}\left(e^{2\hat{i}\theta}\right) \right)  \right] + R_n \\
     &= A_n+h_nB^{(1)}_n+R_n,
\end{align*}
where
\begin{align}
   A_n&=\frac{1}{2}I-\frac{\delta}{2}T_{n-1}\left(e^{\hat{i} \theta}\right)-\frac{\eta}{2}T_{n-1}\left( e^{2 \hat{i} \theta} \right) \\
    B_n^{(1)} &= {\rm diag} \left(-\frac{\psi^{''} \left( \frac{(i-1)T}{n} \right)}{\psi^{'} \left( \frac{(i-1)T}{n} \right)}\right)\left(\frac{1}{4}I+\frac{\eta}{4}T_{n-1}\left(e^{2\hat{i}\theta}\right)\right),
\end{align}
and $\left\{ {R_n \over h_n}\right\}_n$ is the normalized zero-distributed residual.
We observe that the expansion we derived leads to the definition of a momentary GLT expansion since the sequences $\left\{A_n\right\}_n$ and $\left\{B^{(1)}_n\right\}_n$ are GLT. This is due to the fact that the matrices are composed of sums of products of diagonal samplings and the Toeplitz matrices that we involve, so we can apply {\rm[{\bf GLT5}]}. Moreover, setting $h=\frac{1}{n}=\frac{h_n}{T}$, observing that $\hat{\psi}'(x)=T \psi'(xT)$, we infer
\begin{align}
    \left\{ A \right\}_n & \sim_{{\rm GLT}} \kappa(\theta) \\
    \left\{ {B^{(1)}_n} \right\}_n & \sim_{{\rm GLT
    }} l(x, \theta)= -\frac{ \hat{\psi}^{''}(x)}{\hat{\psi}^{'}(x)}\left[\frac{1}{4}+\frac{\eta}{4}e^{2\hat{i}\theta}\right].
\end{align}
Thus, the first order momentary expansion is
\begin{equation}
    \{\mathbb{L}_n\}_n = \{A_n\}_n + h \{B^{(1)}_n\}_n + \{R_n\}_n
\end{equation}
with momentary symbols
\begin{equation}
    f_n(x,\theta)=k(\theta)+h l(x, \theta).
\end{equation}

Finally, the momentary expansion of $\mathbb{S}_n = \Re(\mathbb{L}_n)$ follows directly by the $*$-algebra structure of GLT matrix-sequences and by the linearity of the momentary GLT symbol expansion.
\end{prf*}

\subsection{The decomposition approach}

In the current section we propose a decomposition in rank $p$ - with $p$ at most equal to $2$ - nonnegative definite matrices of the Hermitian matrix $\mathbb{S}_n$.

The decomposition technique has been proven to be very powerful for deducing that Hermitian matrix-valued operators are linear and positive (refer to the notion of matrix-valued LPOs in \cite{ergo,dyadic-LPO}) and for giving, as a consequence, quite refined spectral localization and spectral distributional results \cite{dyadic-cris0,dyadic-LocX}.
For instance, for the fourth-order boundary value problem
\begin{align}\label{prob14}
	(\alpha(x)v''(x))''=\beta(x) \qquad 0<x<1, \,v(0)=v'(0)=v(1)=v'(1)=0,
\end{align}
we can consider the second-order precision central FD formula of minimal bandwidth. The related approximated discrete equations are
\begin{align}\label{discr1}
	\alpha_{i-1}v_{i-2} -2(\alpha_{i-1}+\alpha_i)v_{i-1}+(\alpha_{i-1}+4\alpha_i+\alpha_{i+1})v_i-2(\alpha_{i+1}+\alpha_i)v_{i+1}+\alpha_{i+1}v_{i+2}= h^4\beta(t)
\end{align}
for $i= 2,3,\ldots , n-1$.  The structure of the resulting matrix $A_n(\alpha)$  is
\begin{align*}\label{LTmatrix}
	\begin{bmatrix}
		(\alpha_3+4\alpha_2+\alpha_1) &-2(\alpha_3+\alpha_2) & \alpha_3 & & &\\
		-2(\alpha_2+\alpha_3)& (\alpha_4+4\alpha_3+\alpha_2) &-2(\alpha_4+\alpha_3) & \alpha_4 & &\\
		\alpha_3& -2(\alpha_3+\alpha_4)& (\alpha_5+4\alpha_4+\alpha_3) &-2(\alpha_5+\alpha_4) & \alpha_5& \\
		& & & & & \\
		&\ddots & \ddots& \ddots& \ddots&\ddots\\
		& & & & & \\
		& \alpha_{n-4}&-2(\alpha_{n-4}+\alpha_{n-3})& (\alpha_{n-2}+4\alpha_{n-3}+\alpha_{n-4}) &-2(\alpha_{n-2}+\alpha_{n-3})&\alpha_{n-3} \\
		& &\alpha_{n-3}&-2(\alpha_{n-3}+\alpha_{n-2})& (\alpha_{n-1}+4\alpha_{n-2}+\alpha_{n-3}) &-2(\alpha_{n-1}+\alpha_{n-2})  \\
		& & & \alpha_{n-2}&-2(\alpha_{n-2}+\alpha_{n-1})& (\alpha_{n}+4\alpha_{n-1}+\alpha_{n-2})
	\end{bmatrix}
	\end{align*}
and the associated sequence has GLT nature.
If we substitute $\alpha(t)= 1$, then  the GLT sequence of matrices $A_n(\alpha)$ reduces to the sequence $\{A_n(1)\}_n$ of Toeplitz matrices for which
\begin{align*}
	A_n(1)=\begin{bmatrix}
		6 & -4 & 1 &  & & &  \\
		-4	 & 6 & -4 & 1 & & &\\
		1&-4	 & 6 & -4 & 1& &  \\
		&\ddots & \ddots & \ddots & \ddots & \ddots & \\
		& & 1&-4	 & 6 & -4 & 1  \\
		& & &1&-4	 & 6 & -4   \\
		& & & & 1&-4	 & 6   \\
	\end{bmatrix}.
\end{align*}

Let $\rho(\theta)$ be the trigonometric polynomial associated with the finite difference stencil $(1, -4, 6, -4, 1)$, i.e.,
$$\rho(\theta)=(2-2\cos(\theta))^2.$$
We have $\rho(\theta)=(2-2\cos(\theta))^2=\exp(2i\theta)-4\exp(i\theta)+6-4\exp(-i\theta)+\exp(-2i\theta)$  in the Fourier variable $\theta$,  $\theta\in [-\pi,\pi]$ and hence $A_n(1)=T_n(\rho)$ in the sense of (\ref{eq:fourier_coefficients})-(\ref{toeplitz}).
On the other hand the matrix-sequence associated to the equations (\ref{discr1}) is of GLT type and its GLT symbol is $\alpha(x)\rho(\theta)$, whenever the weight function  $\alpha:[0,1]\to  \mathbb{R}$ is Riemann integrable.
The matrix $A_n(\alpha)$ can be written as $\sum_{j=1}^n \alpha_j D_j$ where $D_j$ have exactly rank $1$ and each $D_j$ is related to the same (rank $1$) stencil of the form $(1, \, -2, \, 1)^T (1, \, -2, \, 1)$, the structure being inherited by the operator in divergence form.
The decomposition allows to prove in a very elementary way that $A_n(\cdot)$ is a matrix-valued LPO. As a consequence, for a positive function $\alpha$, the eigenvalues of $A_n(\alpha)$ lie in the interval $(0, 16\max \alpha)$ with minimal eigenvalue going to zero as $n^{-4}$ with constant depending on the minimum of $\alpha$ and other parameters. We emphasize that such simple findings are not easy to prove by using standard inclusion/exclusion results like the three Gerschgorin Theorems (see \cite{varga} and references therein for an exhaustive account on such kind of results). Furthermore, it is quite direct to prove that
\[
\{A_n^{-1}(b)A_n(a)\}_n\sim_\lambda (a/b, [0,1])
\]
for positive Riemann integrable functions $a,b$ \cite{ergo,dyadic-LocX}.  For related results regarding matrix-value LPOs and the beautiful Korovkin theory see \cite{koro,koro-nambu}.

Here the analysis takes inspiration from the above techniques, but the problem is intrinsically more complicated since the parameters $r_j$, $j=2,\ldots,n$, do not show up in a linear way in the matrix $\mathbb{S}_n$. However, in the following derivations we deduce a parametric decomposition in low rank matrices (of rank at most 2), leading to conditions for checking the positive definiteness of $\mathbb{S}_n$.

More in detail we study the positivity, in the sense of that the related Rayleigh quotient, of the matrix \(2\mathbb{S}_n=2\Re(\mathbb{L}_n)\) through a suitable sum decomposition. We employ an approach analogous to the positive dyadic sum decomposition.

Consider the following class of low-rank diagonal blocks
\[
A_i =
\begin{bmatrix}
0 & \cdots & 0 & 0 & \cdots & 0 \\
& \ddots & \vdots & & \vdots & \\
& & \frac{a_i}{1+r_{i+1}} & -\delta\frac{\sqrt{r_{i+2}}}{1+r_{i+2}} & \\
& & -\delta\frac{\sqrt{r_{i+2}}}{1+r_{i+2}} & \frac{b_i}{1+r_{i+2}} & \\
& \vdots & & & \ddots & \\
0 & \cdots & 0 & 0 & \cdots & 0 \\
\end{bmatrix},
\]
where the \(2 \times 2\) nonzero diagonal block has matrix coordinates \(A_{(i,i+1),(i,i+1)}\), \(i=1,\dots,n-2\); we denote these small blocks by
\(m_i\) for \(i=1,\dots,n-1\). Similarly, consider

\[
B_i =
\begin{bmatrix}
0 & \cdots & 0 & 0 & 0 & \cdots & 0 & 0 \\
& \ddots & \vdots & & \vdots & \vdots & \\
& & \frac{c_i}{1+r_{i+1}} & 0 & -\frac{\eta{\sqrt{r_{i+2}r_{i+3}}}}{1+r_{i+3}} & \\
& & 0 & 0 & 0 & \\
& & -\frac{\eta{\sqrt{r_{i+2}r_{i+3}}}}{1+r_{i+3}} & 0 & \frac{d_i}{1+r_{i+3}} & \\
& \vdots & & & \ddots & \ddots & \\
0 & \cdots & 0 & 0 & \cdots & 0 & 0 \\
\end{bmatrix},
\]
for \(i=1,\dots, n-3\). We denote the nonzero $3\times 3$ block by \(p_i\).

We aim to express:
\begin{equation}\label{at at most 2 dec}
2\Re(\mathbb{L}_n) = \sum_{i=1}^{n-2}A_i + \sum_{i=1}^{n-3}B_i.
\end{equation}
This block decomposition is not a dyadic decomposition in general, since the rank of the small blocks is at most equal to 2. Our objectives are to find the varying coefficients such that:
\begin{enumerate}
    \item The decomposition holds.
    \item The decomposition provides a positive representation of the matrix.
\end{enumerate}

Note that the summands contain symmetric \(2 \times 2\) blocks or \(3 \times 3\) symmetric blocks with the second row and second column being zero and both types of matrices have at most rank $2$ and at least least one eigenvalue positive due to their trace: if we force the nonnegative character of the determinant, then the decomposition will be composed by Hermitian nonnegative definite rank $p$ matrices with $p\in \{1,2\}$.

To verify whether the decomposition holds, we compare the entries of the matrices. A necessary and sufficient condition for equality is given by the following set of constraints
\begin{equation}
    \begin{cases}
        a_1 + c_1 = 2, \\
        a_2 + b_1 + c_2 = 2, \\
        a_{i+1} + b_i + c_{i+1} + d_{i-1} = 2, & \text{for } 2 \leq i \leq n-3, \\
        a_{n-1} + b_{n-2} + d_{n-3} = 2, \\
        b_{n-1} + d_{n-2} = 2.
    \end{cases}
\end{equation}

We now examine the conditions for the positivity of the matrix. A sufficient condition is that all summands are positive semidefinite. To ensure this, we first require that the determinants of \(m_i\) and the \(2 \times 2\) minor of \(p_i\) with nonzero row and column are non-negative that is
\[
\text{d}(m_i) = \frac{a_i b_i}{1+r_{i+1}} - \delta^2\frac{r_{i+2}}{1+r_{i+2}} \ge 0,
\]
\[
\text{d}(p_i^{(2,2)}) = \frac{c_i d_i}{1+r_{i+1}} - \eta^2 \frac{{r_{i+2}r_{i+3}}}{1+r_{i+3}} \ge 0,
\]
where $d(\cdot)$ is the determinant up to positive multiplicative constants which are useless for determinining the global sign.
Furthermore, we require
\[
a_i > 0, \quad c_i > 0.
\]

Note that this is just a typical problem of feasibility of a semi-algebraic set. Considering the vector $x \in \mathbf{R}^{4(n-1)}$
Our condition is equivalent to the following problem: find (if it exists) an $x$ such that
\begin{align}
Mx&=\textbf{2}, \\
    \frac{x_i x_{i+(n-1)+1}}{1+r_{i+1}} - \delta^2\frac{r_{i+2}}{1+r_{i+2}}&\ge0, \quad 1 \le i \le n-2,\\
    \frac{x_{i+{2(n-1)}} x_{i+{3(n-1)+2}}}{1+r_{i+1}} - \eta^2 \frac{{r_{i+2}r_{i+3}}}{1+r_{i+3}} &\ge 0, \quad 1 \le i \le n-3, \\
    x_i &> 0, \quad 1 \le i \le n-2, \\
    x_{i+2(n-1)} &> 0, \quad 1 \le i \le n-3,
\end{align}
with $M$ an underdetermined matrix $M \in \mathbf{R}^{(n-1) \times 4(n-1)}$
\begin{equation*}
    M=\begin{bmatrix}
        M_1 & M_2 & M_3 & M_4
    \end{bmatrix},
\end{equation*}
where the $(n-1) \times (n-1)$ blocks are
\begin{align}
    M_1 &=\text{diag}(1,\dots,1,0), \quad M_2=\text{diag}(0,1,\dots,1), \\
    M_3 &=\text{diag}(1,\dots,1,0,0), \quad M_4=\text{diag}(0,0,1,\dots,1).
\end{align}
In practical approaches, in general this problem can be solved by formulating a semidefinite programming relaxation problem. Once  a feasible $x$ is found, the parameters used for the positive representation are
\begin{align}
    a_i&=x_i, \quad \text{for} \ 1 \le i \le n-2, \\
    b_i&=x_{i+(n-1)+1}, \quad \text{for} \ 1 \le i \le n-2, \\
    c_i&=x_{i+2(n-1)}, \quad \text{for} \ 1\le i \le n-3, \\
    d_i&= x_{i+3(n-1)+2}, \quad \text{for} \ 1 \le i \le n-3.
\end{align}

An algorithm for solving the above problem of feasibility of a semi-algebraic set will give a precise characterization of the positive definiteness of \(2\mathbb{S}_n=2\Re(\mathbb{L}_n)\), with a potential improvement with respect to the analysis in \cite{BDF-Akrivis}.

The numerical tests in fact inform that the conditions of positive definiteness in \cite{BDF-Akrivis} can be improved. From Figure \ref{fig:spectrum_psi_xsquare} we see that $\Re(\mathbb{L}_n)$ is positive definite, but it is clear that some of the $r_i$ (the initial ones) are outside the range given in \cite{BDF-Akrivis}, i.e. 
$r_i \leq 1.9398$, $i=2,\ldots,n$.

\subsection{The limit Toeplitz case}

In this short subsection, we show that the limit Toeplitz case when the grid is equispaced shows also that there is room for improvement with regard to the study in \cite{BDF-Akrivis}. Indeed, we show that in that setting, the generating function of $\mathbb{S}_n$ has positive minimum and hence all its eigenvalues are strictly larger than this minimum and the smallest eigenvalues converge monotonically to it as $n$ tends to infinity, according to the third and fourth parts of Theorem \ref{th:distri-extre}.

More precisely, given the symbol
\[
\Re(\kappa(\theta))=\frac{1}{2}(1-\delta \cos({\theta})-\eta \cos(2\theta)),
\]
our goal is to analyze the extrema of the generating function above. Using the well-known identity $\cos(2 \theta)=2\cos^2(\theta)-1$, we obtain

\begin{equation}
    \Re(\kappa(\theta))=-\eta \cos^2(\theta) - \frac{\delta}{2}\cos(\theta) + \frac{\eta+1}{2}
\end{equation}

By treating $z=\cos(\theta)$ as a variable, we can examine the following general quadratic function

\begin{equation}
P(z)=-\eta z^2 - \frac{\delta}{2}z + \frac{\eta+1}{2}.
\end{equation}

\begin{figure}[!h]
    \centering
    \includegraphics[width=.48\textwidth]{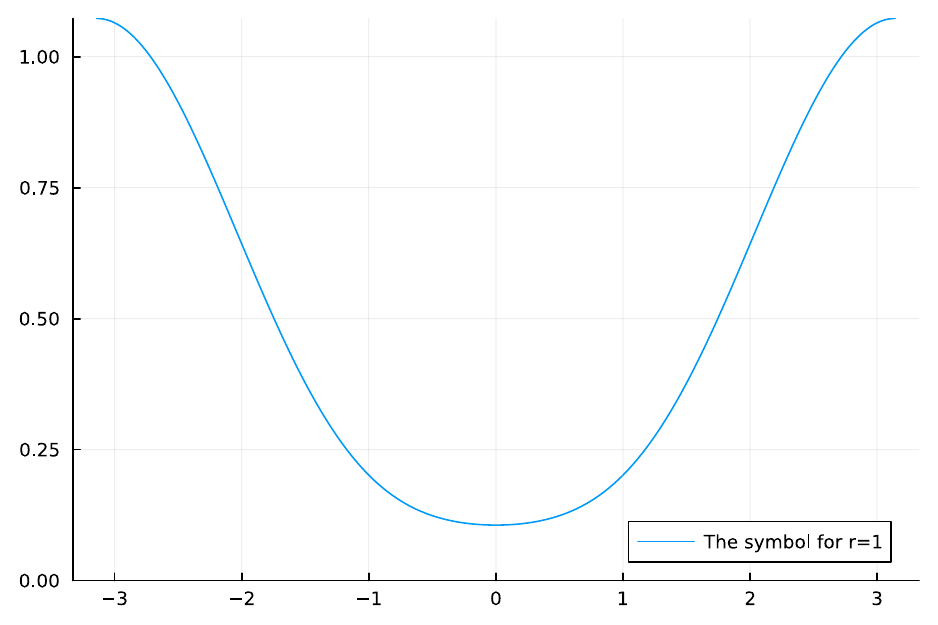}
    \caption{The graph of the generating function (and GLT symbol) of $\mathbb{S}_n(r)=\Re(\mathbb{L}_n(r))=T_{n-1}(\Re(\kappa_r(\theta)))$ for $r=1$. The smoothing parameters are $\delta=0.9672$, $\eta=-0.1793$. The symbol is strictly positive with a minimum value $0.10605$ at zero.}
    \label{fig:symbol_for_r_1}
\end{figure}

\begin{figure}[!h]
    \centering
    \includegraphics[width=.48\textwidth]{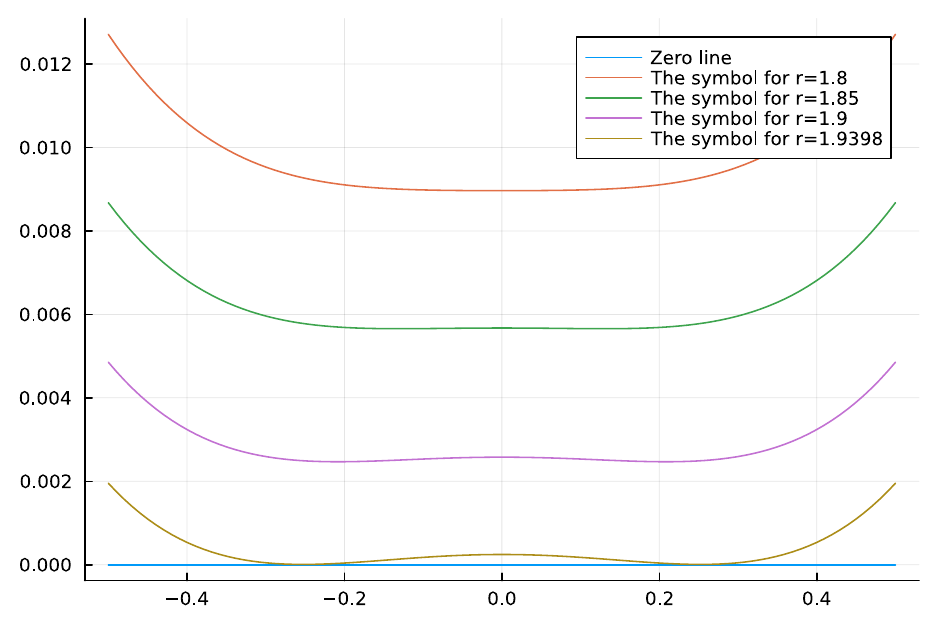}
    \caption{Generating functions (and GLT symbols) of $\mathbb{S}_n(r)=\Re(\mathbb{L}_n(r))=T_{n-1}(\Re(\kappa_r(\theta)))$ for various values of $r$: the largest is the limit value in \cite{BDF-Akrivis} i.e. $r=1.9398$}
    \label{fig:symbol_for_various_r_close}
\end{figure}

\begin{figure}[!h]
    \centering
    \includegraphics[width=.48\textwidth]{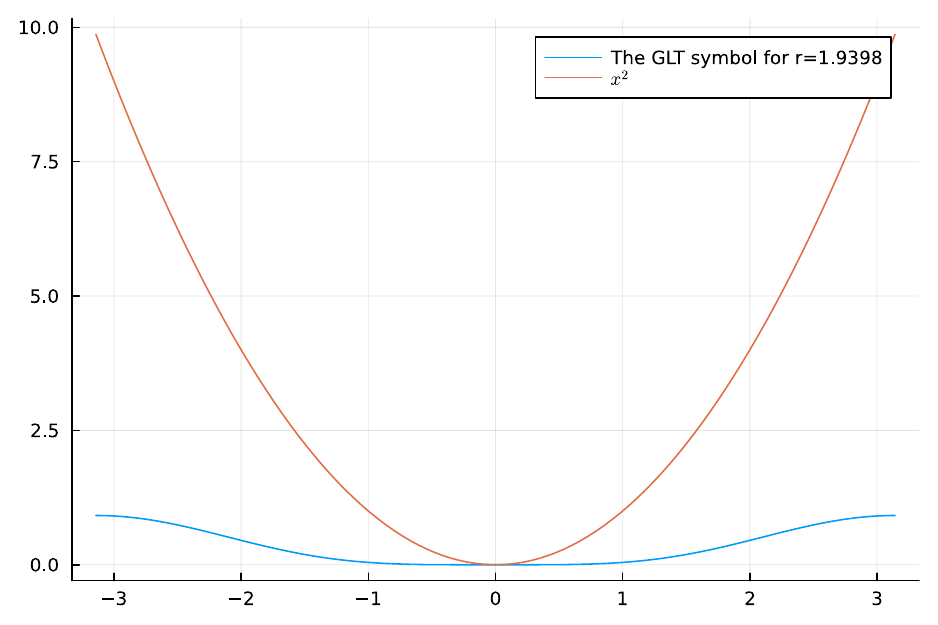}
    \caption{The graph of the generating function (and GLT symbol) of $\mathbb{S}_n(r)=\Re(\mathbb{L}_n(r))=T_{n-1}(\Re(\kappa_r(\theta)))$ for $r=1.9398$ along with the graph of $\theta^2$. The smoothing parameters are $\delta=0.9672$, $\eta=-0.1793$. }
    \label{fig:symbol_vs_xsquare_r19398_distance}
\end{figure}

\begin{figure}[!h]
    \centering
    \includegraphics[width=.48\textwidth]{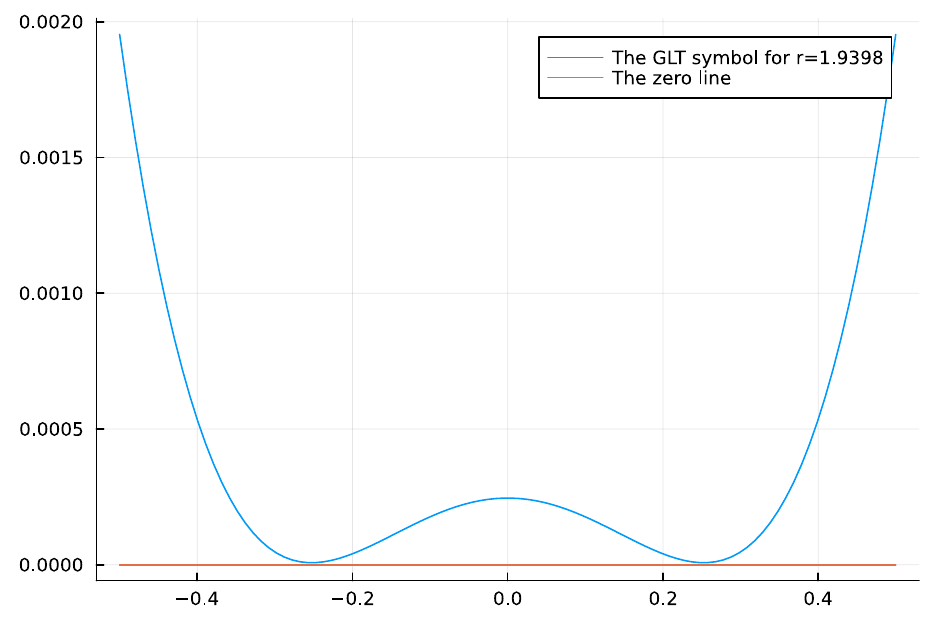}
    \caption{A closer look at the generating function (and GLT symbol) with same parameters as in Figure \ref{fig:symbol_vs_xsquare_r19398_distance} reveals its two minima symmetrically located with respect to $\theta=0$.}
    \label{fig:symbol_vs_zeroline_r19398_close}
\end{figure}

Trivially, the unique extremum \((z_e, p_e)\) of the quadratic function \(y(z) = Az^2 + Bz + C\) is located at the vertex of the corresponding conic section at the point
\[
z_e = -\frac{B}{2A}.
\]
Within our algebraic framework, the extremal point is
\[
z_e= -\frac{\delta}{4\eta}
\]
with extremal value
\begin{equation*}
    P(z_e)=\frac{\delta^2}{16\eta}+(\eta+1)/2.
\end{equation*}
In our particular setting, we are forced to consider the interval $[-1,1]$, where $z=\cos(\theta)$ belongs, with $\mu = -0.1793$ and $b = 0.9672$. For these given parameters, the polynomial is convex and takes the absolute minimum at
\[
x_e = -\frac{0.9672}{4 \times -0.1793} \approx 1.3486,
\]
with minimum value being
\[
P(z_e) \approx 0.0843.
\]
Hence, the global minimum of the polynomial is outside our range of interest, as we are considering a polynomial of cosines, and is also positive. Consequently, we can see that the initial trigonometric symbol is positive.
Precisely, the true global maximum is found at $z_M=-1$ with a maximum value of
\begin{equation}
P(z_M)=1.07325,
\end{equation} while the actual global minimum occurs at $z_m=1$ with a minimum value of
\begin{equation}
P(z_m)=0.10605,
\end{equation}
as it can be seen in Figure \ref{fig:symbol_for_r_1}. Hence the eigenvalues of $\mathbb{S}_n=\Re(\mathbb{L}_n)=T_{n-1}(\Re(\kappa(\theta)))$ are strictly larger than $0.10605$ and the smallest eigenvalues converge monotonically to it as $n$ tends to infinity, according to the third and fourth parts of Theorem \ref{th:distri-extre}. Since in \cite{BDF-Akrivis} the condition of positive definiteness is given as $r_i\le 1.9398$, $i=2,\ldots,n$, here we check the condition under the stationary assumption (equispaced grid) that $r_i=r\le 1.9398$, $i=2,\ldots,n$. For $r$ tending to the limit value $1.9398$, the minimum of the function in Theorem \ref{th:gen} with constant $\phi(x)=r$ reported below
\[
\Re(\kappa_r(\theta))=\frac{1}{r+1}(1-\delta \sqrt{r}\cos({\theta})-\eta r\cos(2\theta)),
\]
moves to zero (see Figure \ref{fig:symbol_for_various_r_close}, as it can be easily checked analytically. Hence, by invoking again the third and fourth parts of Theorem \ref{th:distri-extre}, we deduce that $\mathbb{S}_n(r)=\Re(\mathbb{L}_n(r))=T_{n-1}(\Re(\kappa_r(\theta)))$ is positive definite for every $n$ if and only if $r\le 1.9398$. Conversely, for $r>1.9398$, the minimum of $\Re(\kappa_r(\theta))$ is negative and hence Theorem \ref{th:distri-extre} and Theorem 2.5 in \cite{Sh} imply that
\begin{equation}\label{Sh-conseq-th2.5}
{\rm cardinality}(j: \lambda_j(\mathbb{S}_n(r))<0)=
\frac{n}{\pi}  \mu(\theta\in [0,\pi]:\ \Re(\kappa_r(\theta))<0)
+ O(1)
\end{equation}
with $\mu(\theta:\ \Re(\kappa_r(\theta))<0)>0$ because the cosine polynomial $\Re(\kappa_r(\theta))$ is continuous and has a negative minimum. Figures \ref{fig:symbol_vs_xsquare_r19398_distance} and \ref{fig:symbol_vs_zeroline_r19398_close} give more details on the limit case of $r=1.9398$
with equispaced grid.

As a consequence of the previous study, in particular in the light of (\ref{Sh-conseq-th2.5}), the condition in \cite{BDF-Akrivis} i.e. $r_i=r\le 1.9398$, $i=2,\ldots,n$, is necessary and sufficient for the positive definite character of $\mathbb{S}_n(r)$ for any $n$ in the case of equispaced grids.

In the subsequent part regarding numerical tests, we also report an example with $t_i=\psi(\frac{i}{n})$, $i=0,\ldots,n$, the condition $r_i\le 1.9398$, $i=2,\ldots,n$ is violated but the matrices $\mathbb{S}_n=\Re(\mathbb{L}_n)$ in Theorem \ref{th:specific} are positive definite for every $n$.
%

\subsection{Numerical evidences}

The current subsection containing the numerical tests is divided into two parts: first, we give a set of numerical evidences showing the distributional results in clear way, both in the sense of the eigenvalues and singular values (refer to Theorem \ref{th:gen} and Theorem \ref{th:specific}); then we show numerical tests regarding the positive definiteness of $\mathbb{S}_n=\Re(\mathbb{L}_n)$ in the case of variable $r_i$, $i=2,\ldots,n$, associated to a nonconstant function $\psi(x)=x^2$ and to random values.

\subsubsection{Numerical tests for the distributional results}
\begin{exl}
For this example we use an increasing smooth function $\phi(x)=x^2$ and the parameters $\delta=1$, $\eta=-0.5$. In Figure \ref{fig:spectrum_phi_square} the singular values and the eigenvalues of the generated matrices for $n\in\{80,120\}$ along with GLT symbol of the matrix-sequence are shown, everything rearranged in increasing order. As stated in Theorem \ref{th:gen} the singular values of the matrix-sequence and the eigenvalues of the symmetrized matrix-sequence are distributed as their GLT symbols, in accordance with Figure \ref{fig:spectrum_phi_square}.
\begin{figure}[!h]
    \centering
    \includegraphics[width=.48\textwidth]{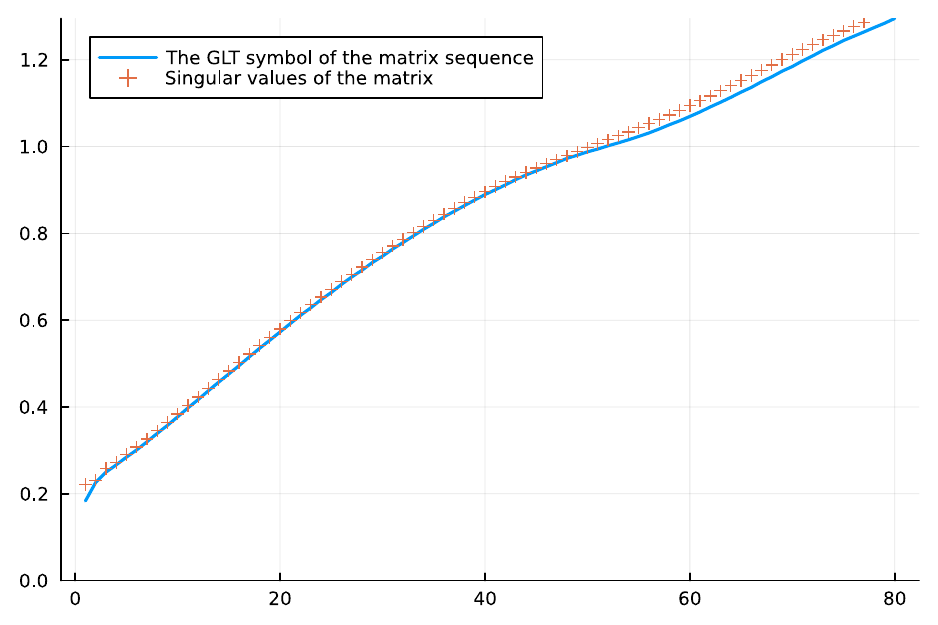}
    \includegraphics[width=.48\textwidth]{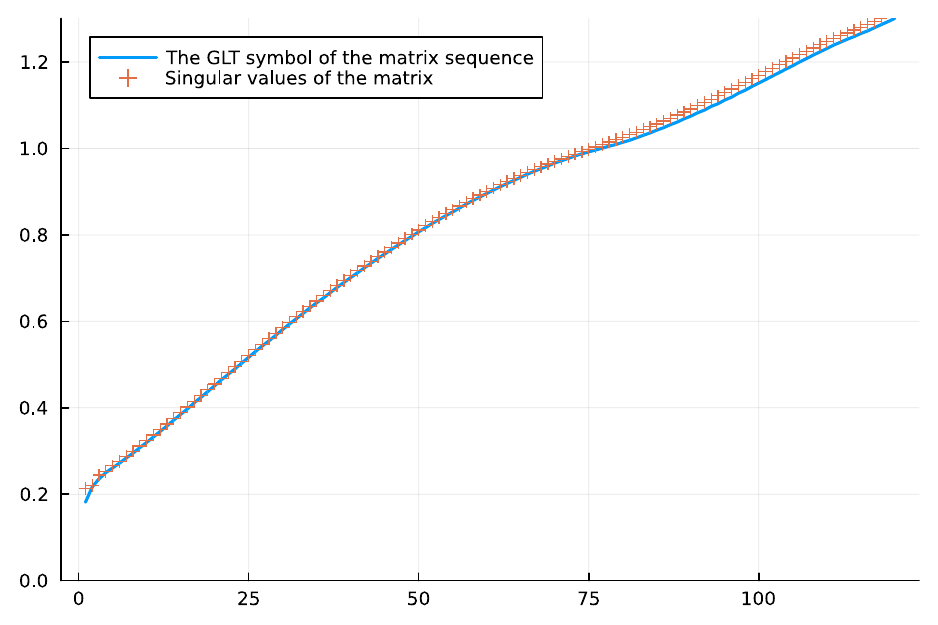}\\
    \includegraphics[width=.48\textwidth]{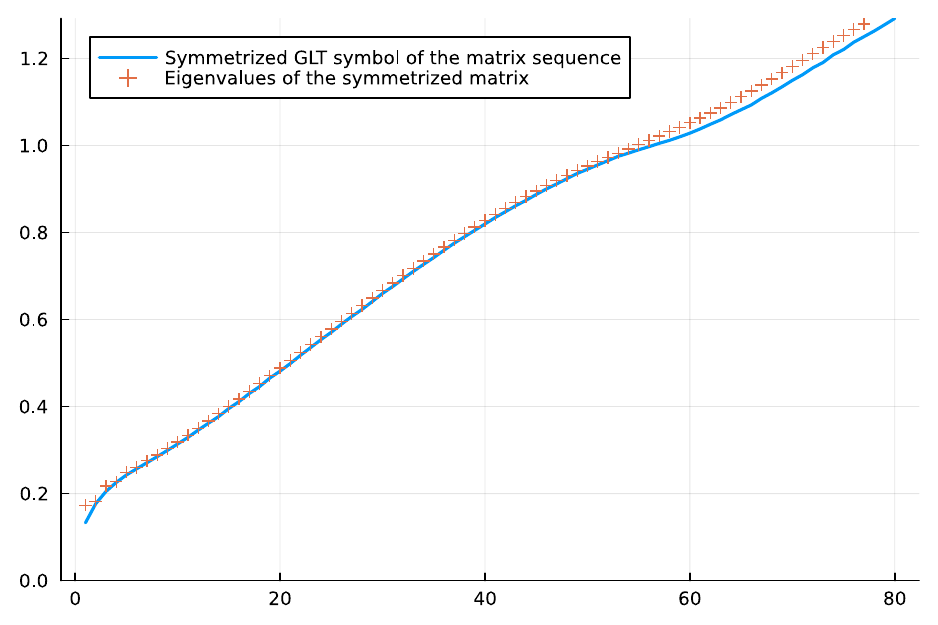}
    \includegraphics[width=.48\textwidth]{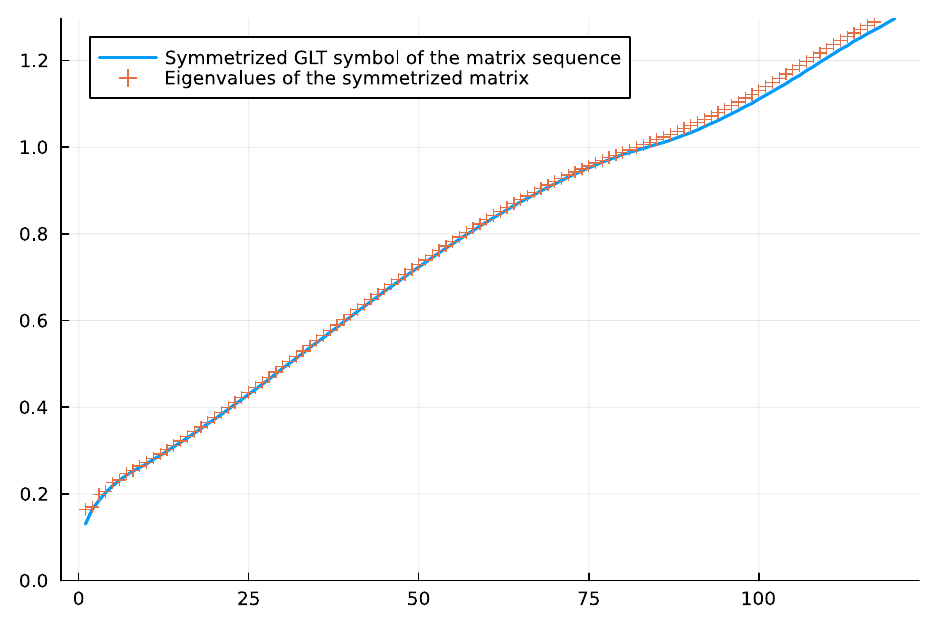}
    \caption{$\phi(x)=x^2$, {$\delta=1$}, {$\eta=-0.5$}
    {Top: The singular values of the matrix for $n=80$ (left) and $n=120$ (right) and the absolute value of the GLT symbol of the matrix-sequence }
    {Bottom: The eigenvalues of the symmetrized matrix for $n=80$ (left) and $n=120$ (right) and the GLT symbol of the matrix-sequence} }
    \label{fig:spectrum_phi_square}
\end{figure}
\end{exl}

\begin{exl}
For this example, we use a non-monotone smooth function $\phi(x) = 1 + \cos(2x)$ and the parameters $\delta = 1$ and $\eta = 1$. In Figure \ref{fig:spectrum_phi_1pluscos2x}, the singular values and the eigenvalues of the generated matrices for $n \in \{80, 120\}$, along with the GLT symbol of the matrix-sequence, are presented, with everything arranged exactly as in the previous example. Again, the singular values of the matrix-sequence and the eigenvalues of the symmetrized matrix-sequence are distributed as their GLT symbols, as predicted in Theorem \ref{th:gen}.
\begin{figure}[!h]
    \centering
    \includegraphics[width=.48\textwidth]{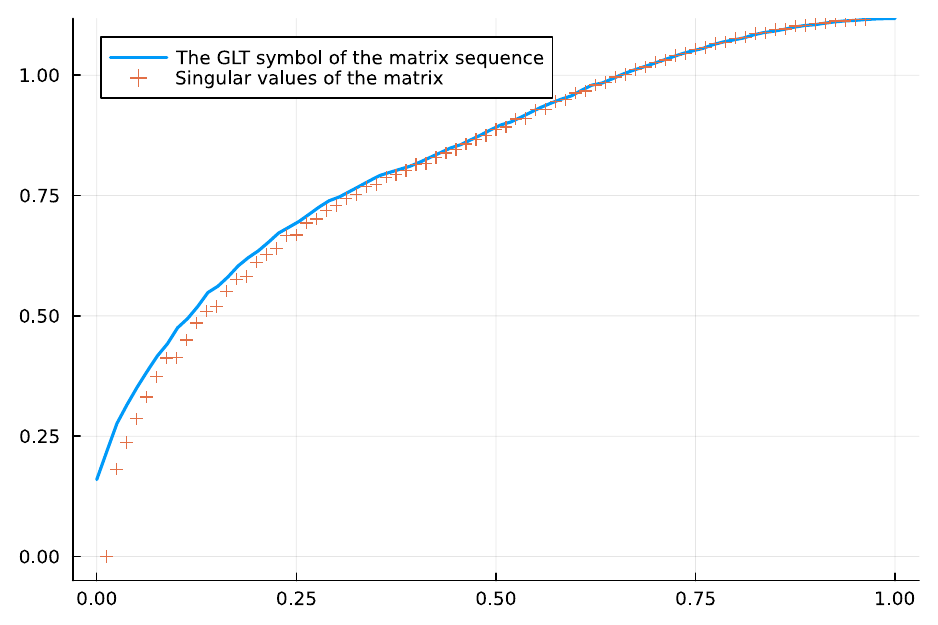}
    \includegraphics[width=.48\textwidth]{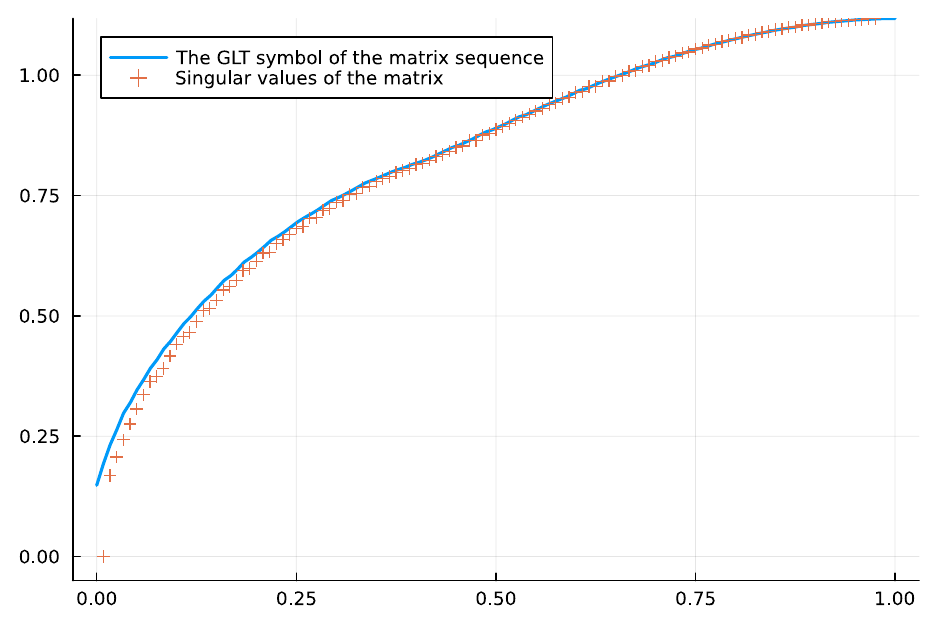}\\
    \includegraphics[width=.48\textwidth]{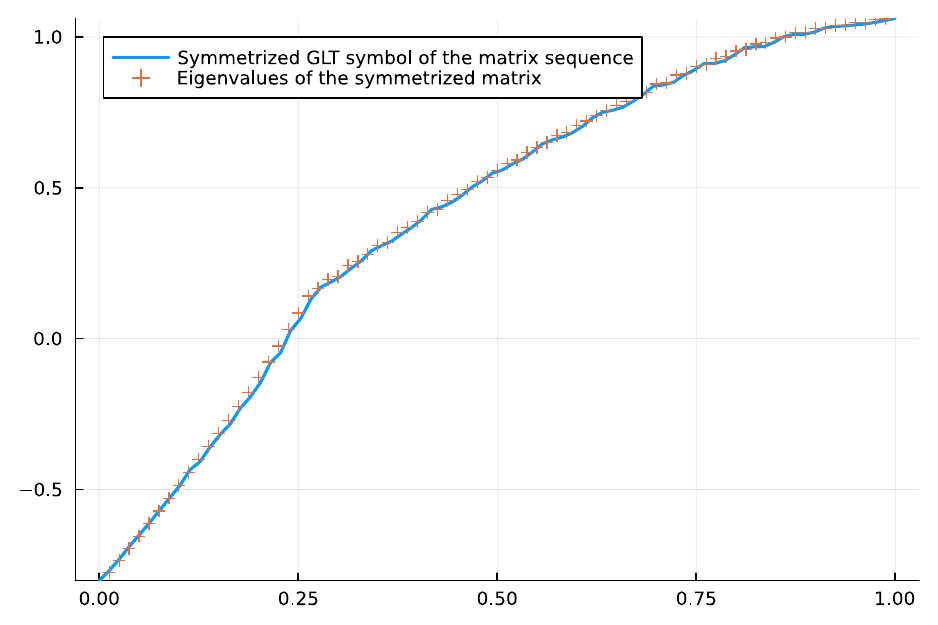}
    \includegraphics[width=.48\textwidth]{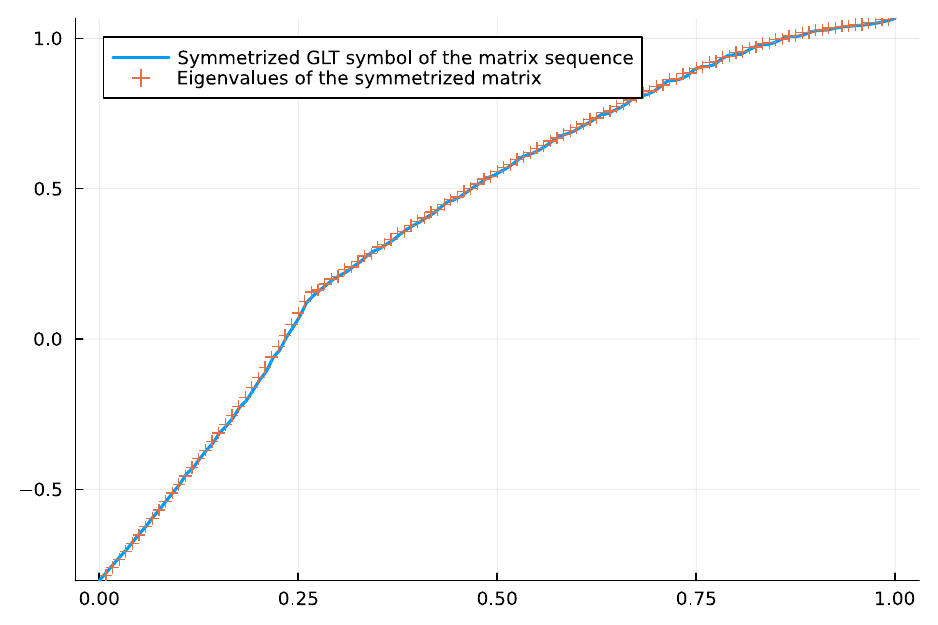}
    \caption{$\phi(x)=1+\cos(2x)$, {$\delta=1$}, {$\eta=1$}
    {Top: The singular values of the matrix for $n=80$ (left) and $n=120$ (right) and the absolute value of the GLT symbol of the matrix-sequence }
    {Bottom: The eigenvalues of the symmetrized matrix for $n=80$ (left) and $n=120$ (right) and the GLT symbol of the matrix-sequence} }
    \label{fig:spectrum_phi_1pluscos2x}
\end{figure}
\end{exl}

\subsubsection{Numerical tests for the nonnegative definiteness}


\begin{exl}
In a more general setting we here consider the case where the time discretization points are such as $t_i=\psi(\frac{i}{n})=(\frac{i}{n})^2$ for $i=0,\dots,n$ and $r_i=\frac{t_{i}-t_{i-1}}{t_{i-1}-t_{i-2}}$. The smoothing parameters are kept $\delta=0.9672$, $\eta=-0.1793$ as suggested in \cite{BDF-Akrivis}. As predicted in (\ref{fata-morgana}) we can see in Figure \ref{fig:spectrum_psi_xsquare} that the case simulates the case where $\phi=1$ (see the previous example). The symmetrized matrices $\mathbb{S}_n=\Re(\mathbb{L}_n)$ remain positive definite for every $n$, although the condition in \cite{BDF-Akrivis} i.e.
\[
r_i \leq 1.9398,
\]
$i=2,\ldots,n$,  is violated.
The smallest eigenvalue gets closer to $0$ because $r_2=\frac{(\frac{2}{n})^2-(\frac{1}{n})^2}{(\frac{1}{n})^2}= 3$ but this is independent of $n$ and the boundary of the smallest eigenvalue is independent of the matrix size. The latter shows numerically that there is room for improving the condition in \cite{BDF-Akrivis}, when $\psi$ is nonconstant and using the decomposition in (\ref{at at most 2 dec}). We finally notice that for a more graded distribution of $t_i=\psi(\frac{i}{n})$ with $\psi(x)=x^3$ the matrices $\mathbb{S}_n=\Re(\mathbb{L}_n)$ become indefinite for $n$ large enough.
\begin{figure}[!h]\label{fig:spectrum_psi_xsquare}
    \centering
    \includegraphics[width=.48\textwidth]{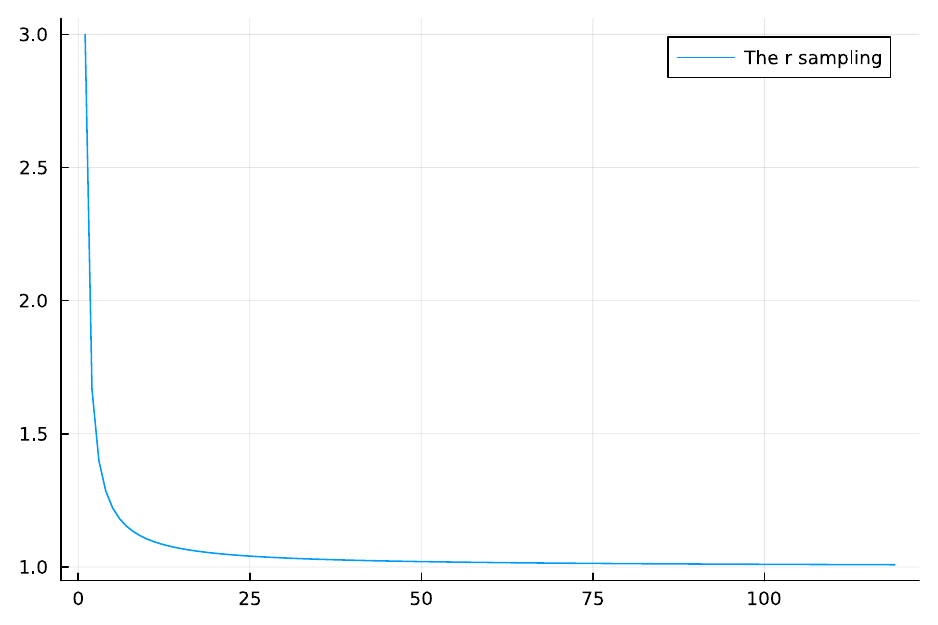}\\
    \includegraphics[width=.48\textwidth]{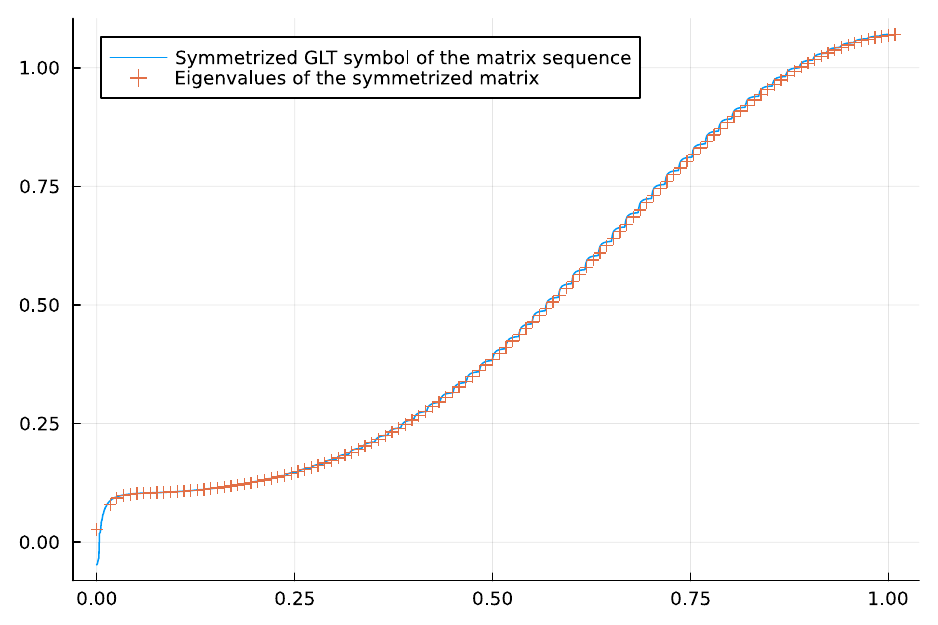}
    \includegraphics[width=.48\textwidth]{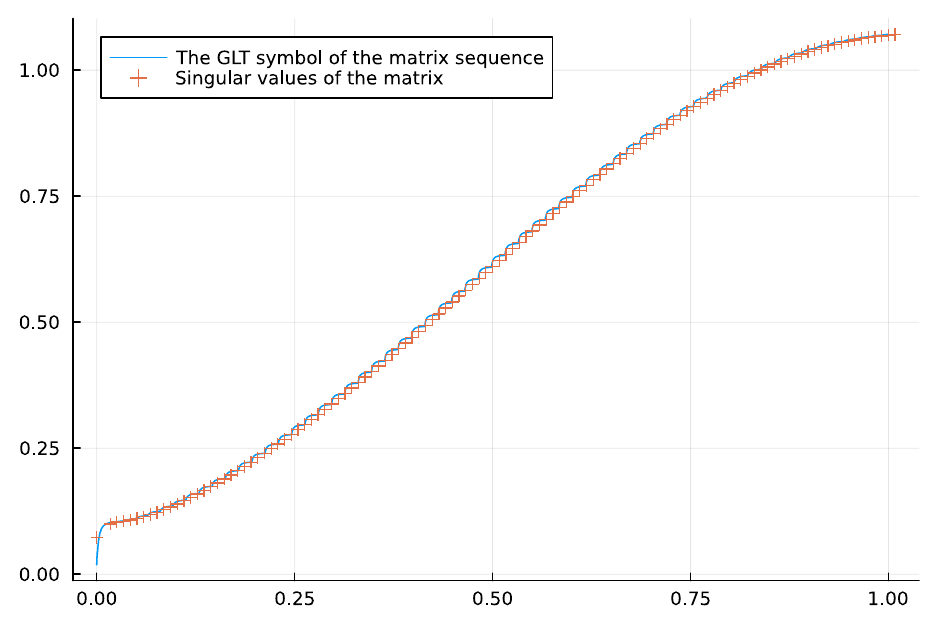}
    \caption{We here consider a time discretization such that $t_i=\psi(\frac{i}{n})=(\frac{i}{n})^2$ for $i=0,\dots,n$ and $r_i=\frac{t_{i}-t_{i-1}}{t_{i-1}-t_{i-2}}$, $i=2,\ldots,n$,  $\delta=0.9672$, $\eta=-0.1793$. The top panel shows the sampling of $r$ while the bottom panels show the eigenvalues (left) and singular values (right) of the generated matrix against the graph of the GLT symbol  }

\end{figure}
\end{exl}

\begin{exl}
In the current example the time discretization points are generated in a random way under the restriction that $0.5\leq \frac{t_{i}-t_{i-1}}{t_{i-1}-t_{i-2}}\leq 1.9398$. The smoothing parameters again are kept $\delta=0.9672$, $\eta=-0.1793$ as suggested in \cite{BDF-Akrivis}. As evident in Figure \ref{fig:spectrum_psi_random} although the values $r_i$, $i=2,\ldots,n$, exhibit a highly erratic behavior,  the symmetrized matrices $\mathbb{S}_n=\Re(\mathbb{L}_n)$ remain positive definite. In this setting, the GLT symbol cannot be defined in a standard way, but only using expected values.
 This aspect goes outside the scope of the present contribution and it is not discussed here. However, a random version of the GLT theory is a field of interest for future researches.
\begin{figure}[!h]
    \centering
    \includegraphics[width=.48\textwidth]{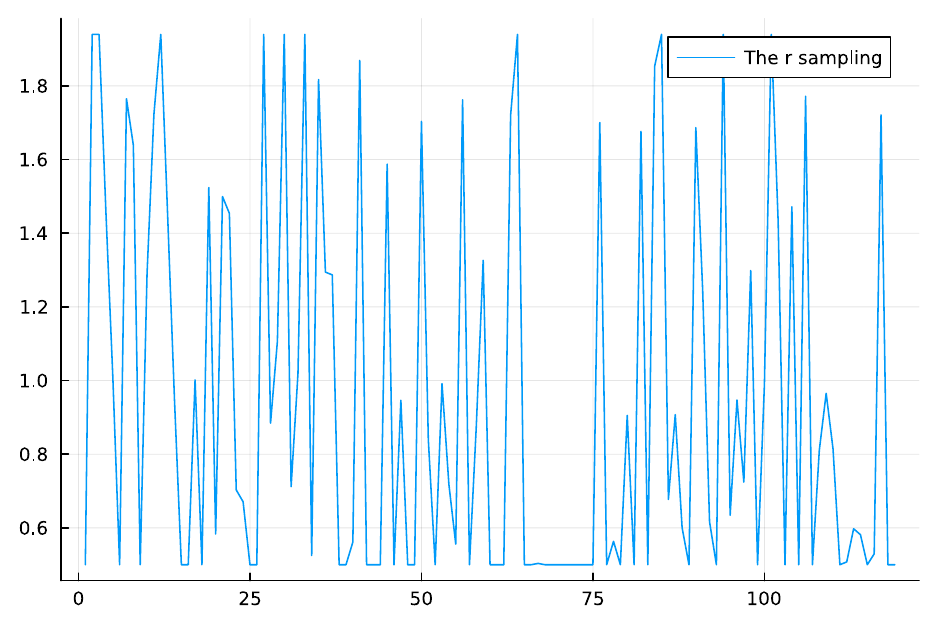}
    \includegraphics[width=.48\textwidth]{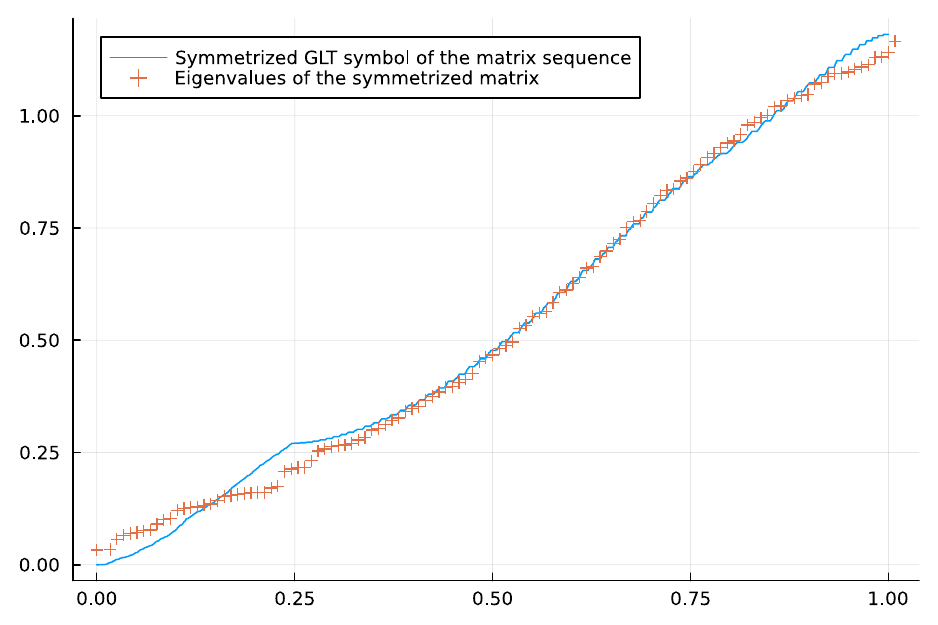}
    \caption{In this last example the time discretization is randomly generated with the restriction that $0.5\leq r_i \leq 1.9398$, $r_i=\frac{t_{i}-t_{i-1}}{t_{i-1}-t_{i-2}}$, $i=2,\ldots,n$, and with $\delta=0.9672$, $\eta=-0.1793$. The symmetrized matrix remains positive definite, with the minimum eigenvalue equal to 0.018985.}
    \label{fig:spectrum_psi_random}
\end{figure}
\end{exl}

\begin{exl}
For giving numerical evidences for the first statement of Theorem \ref{th:specific} we here set $\phi(x)=1$, $\delta=0.9672$, $\eta=-0.1793$. As shown in Figure \ref{fig:spectrum_of_basic_matrix} the  eigenvalues rearranged in increasing order of the symmetrized matrices $\mathbb{S}_n=\Re(\mathbb{L}_n)$ approach the minimum $m_f=0.10605$ and maximum $M_f=1.07375$ of the GLT symbol with a convergence order of 2, as $n$ tends to infinity.

\begin{figure}[!h]
    \centering
    \includegraphics[width=.48\textwidth]{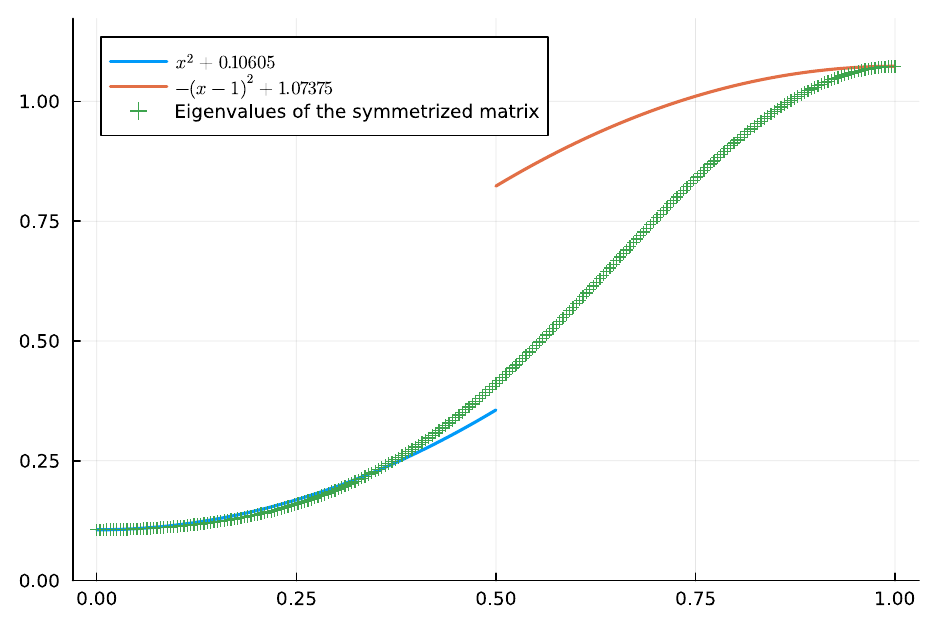}
    \caption{$\phi(x)=1$, $\delta=0.9672$, $\eta=-0.1793$. The smallest eigenvalues of the symmetrized matrices $\mathbb{S}_n=\Re(\mathbb{L}_n)$ approach $m_f=0.10605$ with a convergence order of 2, while the largest eigenvalues similarly approach $M_f=1.07375$ with a convergence order of 2.  }
    \label{fig:spectrum_of_basic_matrix}
\end{figure}
\end{exl}

\section{Conclusions}\label{sec:end}
	
We have considered a general class of matrix-sequences of variable Toeplitz type by proving that they belong to the maximal $*$-algebra of GLT matrix-sequences. We have identified the associated GLT symbols and GLT momentary symbols in the general setting.
Regarding Toeplitz and GLT momentary symbols, new more precise definitions have been given with respect to those in \cite{T-momentary,GLT-momentary}.
When considering variable grid BDF methods for parabolic equations the related structures belong to considered class, and hence using this information we have given a spectral and singular value analysis. Numerical visualizations are also presented corroborating the theoretical analysis.

Finally, we have also proposed a low rank Hermitian nonnegative definite decomposition that, in connection with the notion of matrix-valued LPOs, could lead to precise spectral localization results: this task is however not complete and we leave it for future investigations.

\newpage

\section*{Acknowledgements}
Stefano Serra-Capizzano is partially supported by the Italian Agency INdAM-GNCS. Furthermore, the work of Stefano Serra-Capizzano is funded from the European High-Performance Computing Joint Undertaking  (JU) under grant agreement No 955701. The JU receives support from the European Union’s Horizon 2020 research and innovation programme and Belgium, France, Germany, and Switzerland.
  Stefano Serra-Capizzano is also grateful for the support of the Laboratory of Theory, Economics and Systems – Department of Computer Science at Athens University of Economics and Business.

\clearpage

\newpage

\bibliographystyle{plain}

\end{document}